\documentclass[letter,12pt]{article}%
\usepackage{amsmath,bm}
\usepackage{amsfonts}
\usepackage{amssymb}
\usepackage{graphicx}
\usepackage{rotating}
\usepackage[height=9in,left=1in,right=0.75in,bottom=1in]{geometry}
\usepackage[breaklinks=true,bookmarksopen=true,colorlinks=true,citecolor=blue]%
{hyperref}
\usepackage[authoryear]{natbib}
\usepackage{color}
\usepackage{colortbl}
\usepackage{paralist}
\usepackage{amsmath}%
\providecommand{\U}[1]{\protect\rule{.1in}{.1in}}
\RequirePackage{hypernat}

\catcode`~=11
\newcommand{\urltilde}{\kern -.15em\lower .7ex\hbox{~}\kern .04em}
\catcode`~=13
\def \@seccntformat#1{\csname the#1\endcsname.\quad}
\makeatother
\numberwithin{equation}{section}
\hypersetup{pdftitle={Escanciano}, pdfsubject={Robust Two-Step Semiparametric Estimation}, pdfauthor={Juan Carlos Escanciano}, pdfkeywords={Two-step estimation; Stochastic Equicontinuity;
Semiparametric inference} }
\begin{document}

\title{Cross-Fitting and Fast Remainder Rates for Semiparametric Estimation}
\author{Whitney K. Newey\thanks{Department of Economics, MIT, Cambridge, MA 02139,
U.S.A E-mail: wnewey@mit.edu. }\\\textit{MIT}
\and James M. Robins\thanks{Harvard School of Public Health, 677 Huntington Ave,
Kresge Building Room 823, Boston, MA 02115 U.S.A E-mail:
robins@hsph.harvard.edu. }\\\textit{Harvard}}
\date{January 2018}
\maketitle

\begin{abstract}
There are many interesting and widely used estimators of a functional with
finite semiparametric variance bound that depend on nonparametric estimators
of nuisance functions. We use cross-fitting (i.e. sample splitting) to
construct novel estimators with fast remainder rates. We give cross-fit doubly
robust estimators that use separate subsamples to estimate different nuisance
functions. We obtain general, precise results for regression spline estimation
of average linear functionals of conditional expectations with a finite
semiparametric variance bound. We show that a cross-fit doubly robust spline
regression estimator of the expected conditional covariance is semiparametric
efficient under minimal conditions. Cross-fit doubly robust estimators of
other average linear functionals of a conditional expectation are shown to
have the fastest known remainder rates for the Haar basis or under certain
smoothness conditions. Surprisingly, the cross-fit plug-in estimator also has
nearly the fastest known remainder rate, but the remainder converges to zero
slower than the cross-fit doubly robust estimator. As specific examples we
consider the expected conditional covariance, mean with randomly missing data,
and a weighted average derivative.

\[
\]%
\[
\]

Keywords: Semiparametric estimation, semiparametric efficieny, bias, smoothness.

\end{abstract}

\section{Introduction}

There are many interesting and widely used estimators of a functional with
finite semi-parametric variance bound that depend on the estimation, in a
first step, of nuisance functions, such as conditional expectations or
densities. Examples include estimators of the mean with data missing at
random, the average treatment effect, the expected conditional covariance,
partially linear models, and weighted average derivatives. Because the
nuisance functions can often be high dimensional it is desirable to minimize
the impact of estimating these functions. By using cross-fitting (i.e. sample
splitting) to estimate the nuisance functions we obtain novel estimators whose
second order remainders converge to zero as fast as known possible. In
particular, such estimators are often root-n consistent under minimal
smoothness conditions. Furthermore, such estimators may have higher order mean
square error that converges to zero as fast as known possible.

Bias reduction is key to constructing semiparametric estimators with fast
remainder rates. The rates at which the variance of remainders goes to zero
are quite similar for different semiparametric estimators but the bias rates
differ greatly. We use cross-fitting for bias reduction. We show how fast
remainder rates can be attained by using different parts of an i.i.d. sample
to estimate different components of an estimator.

In this paper we consider regression spline estimation of average linear
functionals of conditional expectations with a finite semiparametric variance
bound, as we have been able to obtain general, precise results for functionals
in this class. The class includes the five examples mentioned above. 

We define a cross fit (CR)$\ $plug-in estimator to be one where we estimate
the functional by simply replacing the unknown conditional expectation by a
nonparametric estimator from a separate part of the sample. Cross-fitting
eliminates an "own observation" bias term, thereby decreasing the size of the
remainder. Functionals in our class have doubly robust influence functions
that depend on two unknown functions. This implies there exists an estimator
depending on both unknown functions that has exact bias zero if the unknown
functions are replaced by fixed functions, at least one of which is equal to
the truth. Here we use double cross-fitting where the two unknown functions
are themselves estimated from separate subsamples, so that the final estimator
depends on three separate subsamples. Surprisingly, single cross fitting in
which both unknown functions are estimated from the same subsample has a
remainder that can converge even slower than CF plug-in estimators. In
contrast, doubly robust estimators with double cross fitting improve on
cross-fit plug-in estimators in the sense that remainder terms can converge at
faster rates. We also show how multiple cross-fitting could be used to reduce
bias for any semiparametric estimator that is a polynomial in first step
spline estmators of unknown functions.

We construct cross-fit (CF) plug-in and doubly cross-fit doubly robust (DCDR)
estimators that are semiparametrically efficient under minimal conditions when
the nuisance functions are in a Holder class of order less than or equal to
one. When a nuisance function is Holder of order exceeding one, we propose
DCDR estimators that have remainders that converge no slower and often faster
than the CF plug-in estimator. In the special case of the expected conditional
covariance functional, the DCDR estimator is always semiparametric efficient
under minimal conditions. For other functionals in our class the CF plug-in
and DCDR estimator are semiparametric efficient under minimal conditions,
provided the conditional expectation is Holder of order greater than or equal
to one-half the regressor dimension; furthermore, in this case, the remainder
goes to zero as fast as known possible for both CF plug-in and DCDR
estimators. When the conditional expectation is Holder of order less than or
equal to one-half the regressor dimension but greater than or equal to one,
the remainder for the DCDR has a remainder that converges faster than the CF
plug-in estimator.

In the case where the conditional expectation is Holder of order no less than
one but less than one-half the regressor dimension, we show semiparametric
efficiency under minimal conditions for the expected conditional covariance,
but not for other functionals. The higher order influence function (HOIF)
estimators of Robins et al. (2008, 2017) and Mukherjee, Newey, and Robins
(2017) will be semiparametric efficient under minimal conditions for these
other functionals, including the mean with data missing at random and the
average treatment effect.

CF plug-in estimators have been considered by Bickel (1982) in the context of
adaptive semiparametric efficient estimation, Powell, Stock, and Stoker (1989)
for density weighted average derivatives, and by many others. Kernel and
series CF plug-in estimators of the integrated squared density and certain
other functionals of a density have been shown to be semiparametric efficient
under minimal conditions by Bickel and Ritov (1988), Laurent (1996), Newey,
Hsieh, and Robins (2004), and Gine and Nickl (2008). Our DCDR estimator
appears to be novel as does the fact that a CF plug-in estimator can be
semiparametric efficient under minimal conditions. Ayyagari (2010), Robins et
al. (2013), Kandasamy et. al. (2015), Firpo and Rothe (2016), and Chernozhukov
et al.(2017) have considered doubly robust estimators that eliminate own
observation terms. Double cross-fitting in double robust estimation appears
not to have been analyzed before.

Our results for splines make use of the Rudelson (1999) law of large numbers
for matrices similarly to Belloni et al.(2015). The results for the CF plug-in
estimator for general splines extend those of Ichimura and Newey (2017) to
sample averages of functionals. The double robustness of the influence
function for the functionals we consider is shown in Chernozhukov et
al.(2016), where the doubly robust estimators of Scharfstein, Rotnitzky, and
Robins (1999), Robins, Rotnitzky, and van der Laan (2000), Robins et. al.
(2008), and Firpo and Rothe (2016) are extended to a wide class of average
linear functionals of expectations.

The DCDR estimator for the mean with missing data and average treatment effect
uses a spline approximation to the reciprocal of the propensity score rather
than the reciprocal of a propensity score estimator. The reciprocal of a
propensity score estimator has been used in much of the previous literature on
plug in and doubly robust estimation, including Robins and Rotnitzky (1995),
Rotnitzky and Robins (1995), Hahn (1998), and Hirano, Imbens, and Ridder
(2003). Estimators based on approximating the reciprocal of the propensity
score have been considered by Robins et al. (2007), Athey, Imbens, and Wager
(2017), and recently in independent work by Hirschberg and Wager (2017).

Other approaches to bias reduction for semiparametric estimators have been
proposed. Robins et al.(2008, 2017) and Mukherjee, Newey, and Robins (2017)
develop higher order influence function (HOIF) estimators with smaller bias.
In Section 2 we will discuss the relationship of this paper to HOIF. Cattaneo
and Jansson (2017) propose promising bootstrap confidence intervals for
plug-in kernel estimators that include bias corrections. Also, Cattaneo,
Jansson, and Ma (2017) show that the jackknife can be used to reduce bias of
plug-in series estimators. For the class of functionals in this paper
cross-fitting removes bias so that there is no need for bootstrap or jackknife
bias corrections in order to attain the fastest remainder rates.

In Section 2 we will describe the cross-fitting approach to bias reduction and
show how it relates to HOIF. Section 3 describes the linear functionals and
regression spline estimators we consider. Sections 4, 5, and 6 give results
for the CF plug-in estimator, the DCDR expected conditional covariance
estimator, and DCDR\ estimators of other linear functionals, respectively.

Before explaining the results of this paper it is helpful to be more specific
about our goal. We will consider i.i.d. data $z_{1},...,z_{n}$. We are
interested in an asymptotically linear semiparametric estimator $\hat{\beta}$
satisfying%
\begin{equation}
\sqrt{n}\left(  \hat{\beta}-\beta_{0}\right)  =\frac{1}{\sqrt{n}}\sum
_{i=1}^{n}\psi\left(  z_{i}\right)  +O_{p}\left(  \Delta_{n}\right)
,\Delta_{n}\longrightarrow0, \label{exp}%
\end{equation}
where $\psi\left(  z\right)  $ is the influence function of $\hat{\beta}$ and
$\Delta_{n}$ characterizes the size of the remainder. Our goal is to find
estimators where $\Delta_{n}$ converges to zero at the fastest known rate.

For the integrated squared density, Bickel and Ritov (1988) gave a kernel
based estimator where the rate for $\Delta_{n}$ is fast enough that
$\hat{\beta}$ is semiparametric efficient under minimal conditions.

To motivate our candidate for the optimal rate the remainder can converge to
zero for series estimators of an average linear functional of a conditional
expectation with positive information bound, we consider the series estimator
of the coefficients of a partially linear regression in Donald and Newey
(1994). The model there is $E[y_{i}|a_{i},x_{i}]=a_{i}^{T}\beta_{0}%
+\lambda_{0}(x_{i})$ where $\lambda_{0}(x_{i})$ is an unknown function of an
$r\times1$ vector $x_{i}$. Consider the estimator $\hat{\beta}$ obtained from
regressing $y_{i}$ on $a_{i}$ and a $\ K\times1$ vector $p(x_{i})$ of power
series or regression splines in an i.i.d. sample of size $n$. Assume that the
functions $\lambda_{0}(x)$ and $\alpha_{0}(x)=E[a_{i}|x_{i}=x]$ are each
members of a Holder class of order $s_{\lambda}$ and $s_{\alpha}$
respectively. Define
\[
\Delta_{n}^{\ast}=\sqrt{n}K^{-(s_{\gamma}+s_{\alpha})/r}+K^{-s_{\gamma}%
/r}+K^{-s_{\alpha}/r}+\sqrt{\frac{K}{n}}.
\]
Donald and Newey (1994) showed that under regularity conditions, including
$K/n\longrightarrow0$, equation (\ref{exp}) is satisfied with $\Delta
_{n}=\Delta_{n}^{\ast}.$ Here $\sqrt{n}K^{-(s_{\gamma}+s_{\alpha})/r}$ gives
the rate at which the bias of $\sqrt{n}(\hat{\beta}-\beta_{0})$ goes to zero.
Also, $K^{-s_{\gamma}/r}$ and $K^{-s_{\alpha}/r}$ are stochastic
equicontinuity bias terms, and $\sqrt{K/n}.$ that accounts for stochastic
equicontinuity and degenerate U-statistic variance terms. Furthermore, there
exists $K=K_{n}$ satisfying $K_{n}/n\longrightarrow0$ such that $\Delta
_{n}^{\ast}\longrightarrow0$ if and only if $s_{\gamma}+s_{\alpha}>r/2$.
However the Donald and Newey (1994) result used the fact that the partially
linear model implies $y_{i}-a_{i}^{T}\beta_{0}$ is mean independent of $a_{i}$
given $x_{i}$ and thus is not a locally nonparametric model. A model is said
to be locally nonparametric if, at each law $P$ in the model, the tangent
space is all of $L_{2}\left(  P\right)  .$Henceforth in this paper, we shall
always assume a locally nonparametric model.

Robins et al. (2009) showed that the condition $s_{\gamma}+s_{\alpha}>r/2$ is
necessary and sufficient for the existence of a semparametric efficient
estimator of
\[
\xi_{0}=E\left[  cov\left(  a_{i},y_{i}|x_{i}\right)  \right]  /E\left[
var(a_{i}|x_{i})\right]  ,
\]
Note $\xi_{0}$ is the probability limit of the Donald and Newey (1994)
estimator regardless of whether the partially linear model holds. That is,
$\xi_{0}$ is the coefficient $b$ in the population linear projection of
$y_{i}$ on all functions of the form $a_{i}b+\lambda(x_{i})$. Robins et al.
(2008) proved sufficiency using a higher order influence function estimator of
$\xi_{0}$, which is a U-statistic whose order increases as $\ln\left(
n\right)  .$ In contrast, the aforementioned estimator of Donald and Newey
(1994), although much simpler, is not semiparametric efficient for $\xi_{0}$
in a locally nonparametric model under the minimal condition $s_{\gamma
}+s_{\alpha}>r/2.$ The current paper was thus motivated by the question
whether one could construct a simple efficient estimator of $\xi_{0}$ whose
remainder $\Delta_{n}$ will go to zero as fast as $\Delta_{n}^{\ast},$ the
fastest rate known to be possible. In summary, our goal is to construct
estimators that are much simpler than the HOIF estimators and yet satisfy
equation (\ref{exp}) with $\Delta_{n}=\Delta_{n}^{\ast}.$

\section{Cross-Fitting and Fast Remainder Rates}

To explain how cross-fitting can help achieve fast remainder rates we consider
estimation of the expected conditional covariance%
\[
\beta_{0}=E[Cov(a_{i},y_{i}|x_{i})]=E[a_{i}\left\{  y_{i}-\gamma_{0}%
(x_{i})\right\}  ],
\]
where $\gamma_{0}(x_{i})=E[y_{i}|x_{i}]$. This object is useful in the
estimation of weighted average treatment effects as further explained below.
We assume that the functions $\gamma_{0}(x)$ and $\alpha_{0}(x)=E[a_{i}%
|x_{i}=x]$ are each members of a Holder class of order $s_{\gamma}$ and
$s_{\alpha}$ respectively.

One way to construct an estimator of $\beta_{0}$ is the \textquotedblleft
plug-in\textquotedblright\ method where a nonparametric estimator $\hat
{\gamma}$ is substituted for $\gamma_{0}$ and a sample average for the
expectation to form%
\[
\bar{\beta}=\frac{1}{n}\sum_{i=1}^{n}a_{i}\{y_{i}-\hat{\gamma}(x_{i})\}.
\]
This estimator generally suffers from an "own observation" bias that is of
order $K/\sqrt{n}$ when $\hat{\gamma}$ is a spline regression estimator, which
converges to zero slower than $\Delta_{n}^{\ast}$. This bias can be eliminated
by replacing $\hat{\gamma}(x)$ with an estimator $\hat{\gamma}_{-i}(x)$ that
does not use $z_{i}$ in its construction. The resulting estimator of
$\beta_{0}$ is
\[
\hat{\beta}=\frac{1}{n}\sum_{i=1}^{n}a_{i}\{y_{i}-\hat{\gamma}_{-i}(x_{i})\}.
\]
This estimator is a cross-fit (CF) plug-in estimator in the sense that
$\hat{\gamma}_{-i}$ uses a subsample that does not include $i$. The
cross-fitting eliminates the own observation bias. The remainder rate
$\Delta_{n}$ for $\hat{\beta}$ will be often be faster than for $\bar{\beta}$,
sometimes as fast as $\Delta_{n}^{\ast}$ as explained below. This approach to
eliminating own observation bias when the first step is a density estimator
has been used by Bickel (1982), Bickel and Ritov (1988), Powell, Stock, and
Stoker (1989), Laurent (1996), and others. Here we obtain the novel result
that, for a spline regression first step, a CF plug-in estimator can have the
fastest rate $\Delta_{n}^{\ast}$ even when the usual plug-in estimator does not.

Doubly robust estimators have another source of bias that can also be
eliminated by double cross-fitting. To explain we consider a single cross-fit
doubly robust estimator of the expected conditional covariance. Let
$\hat{\gamma}_{-i}(x)$ and $\hat{\alpha}_{-i}(x)$ be nonparametric estimators
of $\gamma_{0}(x_{i})=E[y_{i}|x_{i}]$ and $\alpha_{0}(x_{i})=E[a_{i}|x_{i}]$
that do not depend on the $i^{th}$ observation. Consider the estimator%
\[
\check{\beta}=\frac{1}{n}\sum_{i=1}^{n}[a_{i}-\hat{\alpha}_{-i}(x_{i}%
)][y_{i}-\hat{\gamma}_{-i}(x_{i})].
\]
This estimator is doubly robust in the sense of Scharfstein, Rotnitzky, and
Robins (1999) and Robins, Rotnitzky, and van der Laan (2000), being consistent
if either $\hat{\alpha}_{-i}$ or $\hat{\gamma}_{-i}$ are consistent. It uses
cross-fitting to eliminate own observation bias$.$ This estimator does have a
nonlinearity bias since $\hat{\alpha}_{-i}(x_{i})$ and $\hat{\gamma}%
_{-i}(x_{i})$ are constructed from the same data in single crossfitting. That
bias is of the same order $K/\sqrt{n}$ as the own observation bias for a
spline regression plug-in estimator. This bias can be thought of as arising
from nonlinearity of $\check{\beta}$ in the two nonparametric estimators
$\hat{\alpha}_{-i}(x_{i})$ and $\hat{\gamma}_{-i}(x_{i}).$

One can remove the nonlinearity bias in the doubly robust estimator by using
different parts of the data to construct the two nonparametric estimators. Let
$\hat{\gamma}_{-i}(x_{i})$ be constructed from a subset of the observations
that does not include observation $i$ and let $\tilde{\alpha}_{-i}(x_{i})$ be
constructed from a subset of the observations that does not include $i$ or any
observations used to form $\hat{\gamma}_{-i}$. A doubly cross-fit doubly
robust estimator (DCDR) is%
\[
\tilde{\beta}=\frac{1}{n}\sum_{i=1}^{n}[a_{i}-\tilde{\alpha}_{-i}%
(x_{i})][y_{i}-\hat{\gamma}_{-i}(x_{i})].
\]
This estimator uses cross-fitting to remove both the own observation and the
nonlinearity biases. We will show that $\Delta_{n}^{\ast}=\Delta_{n}$ when
$\tilde{\alpha}_{-i}(x_{i})$ and $\hat{\gamma}_{-i}(x_{i})$ are spline
regression estimators for a $K\times1$ vector of multivariate splines of at
least order $\max\{s_{\gamma},s_{a}\}-1$ with evenly spaced knots.
Consequently, this estimator will be root-n consistent and semiparametric
efficient when $s_{\gamma}+s_{\alpha}>r/2$ and $K$ is chosen appropriately,
which is the minimal condition of Robins et al. (2009).

Remarkably, the doubly robust estimator $\check{\beta}$ where $\hat{\alpha
}_{-i}(x_{i})$ and $\hat{\gamma}_{-i}(x_{i})$ use the same data may have a
slower remainder rate than the CF plug-in estimator $\hat{\beta}$. The use of
the same data for $\hat{\alpha}_{-i}(x_{i})$ and $\hat{\gamma}_{-i}(x_{i})$
introduces a bias term of size $K/\sqrt{n}$. Such a term is not present in the
CF plug-in estimator. The $K/\sqrt{n}$ term is eliminated for the doubly
robust estimator by forming $\tilde{\alpha}_{-i}(x_{i})$ and $\hat{\gamma
}_{-i}(x_{i})$ from different samples. We find that the DCDR estimator
$\tilde{\beta}$ improves on the CF plug in estimator by increasing the rate at
which a certain part of $\Delta_{n}$ goes to zero. Specifics will be given below.

We note that the own observation bias can also be thought of as nonlinearity
bias. The parameter $\beta_{0}$ has the form%
\[
\beta_{0}=\int a\{y-\gamma_{0}(x)\}F_{0}(dz),
\]
where $F_{0}$ denotes the distribution of $z=(y,a,x).$ This object is
quadratic in $\gamma_{0}$ and $F_{0}$ jointly. The own observation bias can be
thought of as a quadratic bias resulting from using all the data to
simultaneously estimate $\gamma_{0}$ and the distribution $F_{0}$ of a single
observation. The CF plug in estimator $\hat{\beta}$ eliminates this
nonlinearity bias. Also, the doubly robust estimator can be thought of as
estimating $\int[a-\alpha_{0}(x)][y-\gamma_{0}(x)]F_{0}(dz),$ which is cubic
in $\alpha_{0},$ $\gamma_{0}$, and $F_{0}$ jointly. The DCDR estimator can be
thought of as eliminating the cubic bias by estimating each of $\alpha
_{0}(x),$ $\gamma_{0}(x)$, and $F_{0}$ from distinct groups of observations.

One potential concern about DCDR estimators is that each of the nonparametric
components $\hat{\gamma}$ and $\tilde{\alpha}$ only use a fraction of the data
because they are each based on subsamples that the other does not use. For
example, they only use less than half the data if they are based on
approximately the same subsample size. This does not affect remainder rates
but could affect small sample efficiency. One might be able to improve small
sample efficiency by averaging over DCDR estimators that use different sample
splits to construct $\hat{\gamma}$ and $\tilde{\alpha}$, though that is beyond
the scope of this paper. Our concern in this paper is remainder rates for
asymptotically efficient estimation. 

Cross-fitting can be applied to eliminate bias terms for any estimator that
depends on powers of nonparametric estimators. Such cross-fitting would
replace each power by a product of nonparametric estimators that are computed
from distinct subsamples of the data, analogously to the DCDR estimators above.

We now provide a more quantitative version of our results. Let $p(x)$ be a
vector of multivariate regression splines of dimension $K$ with evenly spaced
knots. We will always take $K=K_{n}$ to satisfy $K\ln\left(  K\right)
/n\rightarrow0.$ Suppose that $\hat{\gamma}_{-i}(x)=p(x)^{T}[\Sigma
_{j\in\mathcal{I}_{\ell}}p(x_{j})p(x_{j})^{T}]^{-1}\Sigma_{j\in\mathcal{I}%
_{\ell}}p(x_{j})y_{j}$ is a series estimator from regressing $y_{j}$ on
$p(x_{j})$ in a subsample of observations indexed by $\mathcal{I}_{\ell}$,
where $\left\{  \mathcal{I}_{\ell}\right\}  _{\ell=1}^{L}$ is a partition of
$\{1,...,n\},$ $i\notin\mathcal{I}_{\ell},$ $L$ is fixed and the number of
elements of each $\mathcal{I}_{\ell}$ is of order $n$. Suppose that for the
doubly robust estimator $\tilde{\alpha}(x_{i})$ is constructed analogously
from a separate subsample.

When $s_{\gamma}\leq1$ and $s_{\alpha}\leq1$ and $p(x)$ is a Haar basis of
dummy variables that are indicator functions of cubes partitioning the support
of $x_{i}$ we show that the CF plug-in estimator has $\Delta_{n}=\Delta
_{n}^{\ast}+\ln(n)K^{-s_{\gamma}/r}$ and the DCDR doubly robust estimator has
$\Delta_{n}=\Delta_{n}^{\ast}.$ Hence the DCDR estimator has the fast
remainder rate. Further the CF plug-in estimator has the fast remainder
$\Delta_{n}^{\ast},$ except at those laws where $K^{-s_{\gamma}/r}$ is the
dominating term in $\Delta_{n}^{\ast}$. At such laws, the DCDR estimator
improves on the CF\ plug-in but only by a factor of $\ln(n).$ We also show
that these results extend to the entire class of average linear functionals of
a conditional expectation with finite semiparametric variance bound.

When $s_{\gamma}$ and $s_{\alpha}$ are any positive numbers and $p(x)$ is a
spline basis of order at least $\max\{s_{\gamma},s_{\alpha}\}-1$ we show that
the CF plug in estimator of the expected conditional covariance has
$\Delta_{n}=\Delta_{n}^{\ast}+\sqrt{K\ln(K)/n}K^{1/2-s_{\gamma}/r}$ and the
DCDR estimator has $\Delta_{n}=\Delta_{n}^{\ast}$.Here the plug-in estimator
has the fast remainder $\Delta_{n}=\Delta_{n}^{\ast}$ for $s_{\gamma}>r/2$ and
the doubly robust estimator has $\Delta_{n}=\Delta_{n}^{\ast}$ for all
$s_{\gamma}$. For other functionals in our class we show that the DCDR
estimator has $\Delta_{n}=\Delta_{n}^{\ast}+\sqrt{K^{3}\ln(K)^{2}/n^{3}%
}K^{1/2-s_{\gamma}/r},$ which has $\Delta_{n}=\Delta_{n}^{\ast}$ when
$[K\ln(K)/n]K^{1/2-s_{\gamma}/r}\longrightarrow0\ .$ Thus the DCDR estimator
has remainder that can converge to zero at a faster rate that of the CF
plug-in estimator.

We note that the source of the term in $\Delta_{n}$ that is added to
$\Delta_{n}^{\ast}$ in each case can be attributed to estimators of the second
moment matrix $\Sigma=E[p(x_{i})p(x_{i})^{T}]$ of the regression splines. If
each $\hat{\Sigma}_{\ell}$ were replaced by $\Sigma$ in the estimators then
the resulting objects would all have $\Delta_{n}=\Delta_{n}^{\ast}.$

For brevity, we demonstrate this only for plug-in estimator. Consider the
plug-in object $\dot{\beta}$ having the same formula as $\hat{\beta}$ except
that $\hat{\gamma}_{-i}(x)$ is replaced by $\dot{\gamma}_{-i}(x)=p(x)^{T}%
\Sigma^{-1}\sum_{j\in\mathcal{I}_{\ell}}p(x_{i})y_{i}/n_{\ell}.$ Let
$\bar{\alpha}(x)=p(x)^{T}\Sigma^{-1}E[p(x_{i})\alpha_{0}(x_{i})]$. Standard
approximation properties of splines give the approximation rates
$\{E[\{\gamma_{0}(x_{i})-\bar{\gamma}(x_{i})\}^{2}]\}^{1/2}=O(K^{-s_{\gamma
}/r})$ and $\{E[\{\alpha_{0}(x_{i})-\bar{\alpha}(x_{i})\}^{2}]\}^{1/2}%
=O(K^{-s_{\alpha}/r}).$ By the Cauchy-Schwartz inequality
\begin{align*}
\sqrt{n}E[\{\alpha_{0}(x_{i})-\bar{\alpha}(x_{i})\}\{\gamma_{0}(x_{i}%
)-\bar{\gamma}(x_{i})\}] &  \leq\sqrt{n}\{E[\{\alpha_{0}(x_{i})-\bar{\alpha
}(x_{i})\}^{2}]\}^{1/2}\{E[\{\gamma_{0}(x_{i})-\bar{\gamma}(x_{i}%
)\}^{2}]\}^{1/2}\\
&  =O(\sqrt{n}K^{-(s_{\gamma}+s_{\alpha})/r}).
\end{align*}
Note also that $E[\dot{\gamma}_{-i}(x)]=\bar{\gamma}(x)=p(x)^{T}\Sigma
^{-1}E[p(x_{i})\gamma_{0}(x_{i})].$ Then the root-n normalized bias of
$\dot{\beta}$ is%
\begin{align}
E\left[  \sqrt{n}\left(  \dot{\beta}-\beta_{0}\right)  \right]   &  =\sqrt
{n}\int a\{y-E\left[  \dot{\gamma}_{-i}(x)\right]  \}F_{0}\left(  dz\right)
-E[a_{i}\{y_{i}-\gamma_{0}(x_{i})\}]\nonumber\\
&  =\sqrt{n}E[a_{i}\{\gamma_{0}(x_{i})-\bar{\gamma}(x_{i})\}]=\sqrt{n}%
E[\alpha_{0}(x_{i})\{\gamma_{0}(x_{i})-\bar{\gamma}(x_{i})\}]\label{bias}\\
&  =\sqrt{n}E[\{\alpha_{0}(x_{i})-\bar{\alpha}(x_{i})\}\{\gamma_{0}%
(x_{i})-\bar{\gamma}(x_{i})\}]=O(\sqrt{n}K^{-(s_{\gamma}+s_{\alpha}%
)/r}),\nonumber
\end{align}
which has our desired $\Delta_{n}^{\ast}$ rate. Also, there will be stochastic
equicontinuity bias terms of order $K^{-s_{\gamma}/r}$ and $K^{-s_{\alpha}/r}$
and stochastic equicontinuity variance and degenerate U-statistic variance
terms of order $\sqrt{K/n}$. Overall the remainder for $\dot{\beta}$ will
satisfy $\Delta_{n}=\Delta_{n}^{\ast}$. Thus, a CF plug-in object $\dot{\beta
}$ where $\Sigma$ replaces each $\hat{\Sigma}_{\ell}$ will have the fast
remainder rate.

We note that the bias in equation (\ref{bias}) depends on the product
$K^{-(s_{\gamma}+s_{\alpha})/r}$ of the approximation rate $K^{-s_{\gamma}/r}$
for $\gamma_{0}(x)$ and the approximation rate $K^{-s_{\alpha}/r}$ for
$\alpha_{0}(x),$ rather than just the bias rate $K^{-s_{\gamma}/r}$ for the
nonparametric estimator being plugged-in. This product form results from the
fact that the parameter of interest $\beta_{0}$ has a finite semiparametric
variance bound. The product bias form in equation (\ref{bias}) for plug-in
series estimators was shown in Newey (1994).

It is interesting to compare our estimators with HOIF estimators. We continue
to focus on the average conditional covariance. The HOIF estimator of that
$\beta_{0}$ can depend on initial estimators $\hat{\gamma}(x)$ and
$\hat{\alpha}(x)$ of $\gamma_{0}(x)$ and $\alpha_{0}(x)$ obtained from a
training subsample. For a vector of spline regressors $p(x)$ let $\hat{\Sigma
}$ be the sample second moment matrix of $p(x)$ from the training sample. Let
$\hat{B}(x)=\hat{\Sigma}^{-1}[p(x)p(x)^{T}-\hat{\Sigma}]$ and%
\begin{align*}
\hat{\beta}_{H}  &  =\frac{1}{n}\sum_{i=1}^{n}[a_{i}-\hat{\alpha}%
(x_{i})][y_{i}-\hat{\gamma}(x_{i})]-\frac{1}{n(n-1)}\sum_{i\neq j}[a_{i}%
-\hat{\alpha}(x_{i})]p(x_{i})^{T}\hat{\Sigma}^{-1}p(x_{j})[y_{j}-\hat{\gamma
}(x_{j})]\\
&  +\sum_{q=1}^{Q}\frac{(-1)^{q+1}(n-2-q)!}{n!}\sum_{i\neq j}[a_{i}%
-\hat{\alpha}(x_{i})]p(x_{i})^{T}\left[  \sum_{\ell_{1}\neq\cdots\neq\ell
_{q}\neq i\neq j}\Pi_{r=1}^{q}\hat{B}(x_{\ell_{r}})\right]  \hat{\Sigma}%
^{-1}p(x_{j})[y_{j}-\hat{\gamma}(x_{j})],
\end{align*}
where all the sums are over an estimation subsample that does not overlap with
the training sample. This $\hat{\beta}_{H}$ is the empirical HOIF estimator of
Mukherjee, Newey, and Robins (2017) of order $Q+2$. By Theorem 3 of Mukherjee,
Newey, and Robins (2017) the bias of $\sqrt{n}(\hat{\beta}_{H}-\beta_{0})$
conditional on the training sample has order%
\[
\sqrt{n}\left\Vert \hat{\alpha}-\alpha_{0}\right\Vert _{2}\left\Vert
\hat{\gamma}-\gamma_{0}\right\Vert _{2}\left(  \frac{K\ln(K)}{n}\right)
^{Q/2}=\left\Vert \hat{\alpha}-\alpha_{0}\right\Vert _{2}\left\Vert
\hat{\gamma}-\gamma_{0}\right\Vert _{2}K\ln(K)\left(  \frac{K\ln(K)}%
{n}\right)  ^{(Q-1)/2}.
\]
where $\left\Vert \delta\right\Vert _{2}=\{E[\delta(x_{i})^{2}]\}^{1/2}.$ The
order of this bias will be smaller than $\sqrt{K/n}$ as long as $K$ grows no
faster than $n^{1-\varepsilon\text{ }}$for some $\varepsilon>0$, although that
is not needed for semiparametric efficiency. As shown in Mukherjee, Newey, and
Robins (2017), if $Q$ grows like $\sqrt{\ln(n)},$ $K$ like $n/\ln(n)^{3},$ and
other regularity conditions are satisfied then $\hat{\beta}_{H}$ will be
semiparametric efficient under the minimal condition $s_{\gamma}+s_{\alpha
}>r/2$ of Robins et al.(2009).

We can explain the different properties of HOIF and series estimators by
comparing the CF plug-in estimator with the HOIF when the training sample
estimators $\hat{\gamma}$ and $\hat{\alpha}$ are set equal to zero. In that
case the HOIF estimator is%
\begin{align*}
\hat{\beta}_{H}  &  =\frac{1}{n}\sum_{i=1}^{n}a_{i}y_{i}-\frac{1}{n(n-1)}%
\sum_{i\neq j}a_{i}p(x_{i})^{T}\hat{\Sigma}^{-1}p(x_{j})y_{j}\\
&  +\sum_{q=1}^{Q}\frac{(-1)^{q+1}(n-2-q)!}{n!}\sum_{i\neq j}a_{i}p(x_{i}%
)^{T}\left[  \sum_{\ell_{1}\neq\cdots\neq\ell_{q}\neq i\neq j}\Pi_{r=1}%
^{q}\hat{B}(x_{\ell_{r}})\right]  \hat{\Sigma}^{-1}p(x_{j})y_{j}.
\end{align*}
Consider $\check{\gamma}_{-i}(x)=p(x)^{T}\hat{\Sigma}^{-1}\sum_{j\neq
i}p(x_{j})y_{j}/(n-1).$ This is an estimator of $\gamma_{0}(x)$ that is like a
series estimator except the inverse second moment matrix $\hat{\Sigma}^{-1}$
comes from the training sample and the cross-moments $\sum_{j\neq i}%
p(x_{j})y_{j}/(n-1)$ from the estimation subsample. The first two terms of the
HOIF estimator can then be written as%
\[
\check{\beta}=\frac{1}{n}\sum_{i=1}^{n}a_{i}[y_{i}-\check{\gamma}_{-i}%
(x_{i})].
\]
Let $T$ denote the training sample. Then we have
\begin{align*}
E[\check{\beta}-\beta_{0}|T]  &  =E[\alpha_{0}(x_{i})\{\gamma_{0}%
(x_{i})-\check{\gamma}_{-i}(x_{i})\}]=E[\alpha_{0}(x_{i})\gamma_{0}%
(x_{i})]-E[\alpha_{0}(x_{i})p(x_{i})^{T}]\hat{\Sigma}^{-1}E[p(x_{i})\gamma
_{0}(x_{i})]\\
&  =E[\alpha_{0}(x_{i})\gamma_{0}(x_{i})-\bar{\alpha}(x_{i})\bar{\gamma}%
(x_{i})]+E[\alpha_{0}(x_{i})p(x_{i})^{T}](\Sigma^{-1}-\hat{\Sigma}%
^{-1})E[p(x_{i})\gamma_{0}(x_{i})]\\
&  =O(K^{-(s_{\gamma}+s_{\alpha})/r})+\Lambda(\hat{\Sigma},\Sigma
),\Lambda(\hat{\Sigma},\Sigma)=E[\alpha_{0}(x_{i})p(x_{i})^{T}](\Sigma
^{-1}-\hat{\Sigma}^{-1})E[p(x_{i})\gamma_{0}(x_{i})].
\end{align*}
Thus the bias of $\check{\beta}$ is the sum of the approximation bias
$K^{-(s_{\gamma}+s_{\alpha})/r}$ and $\Lambda(\hat{\Sigma},\Sigma).$ The rest
of the HOIF estimator, i.e. $\hat{\beta}_{H}-\check{\beta}$, can be thought of
as a bias correction for $\Lambda(\hat{\Sigma},\Sigma).$ Note that
\[
E[\hat{\beta}_{H}-\check{\beta}|T]=\sum_{q=1}^{Q}\frac{(-1)^{q+1}(n-2-q)!}%
{n!}E[\alpha_{0}(x_{i})p(x_{i})]^{T}\left[  \hat{\Sigma}^{-1}(\Sigma
-\hat{\Sigma})\right]  ^{q}\hat{\Sigma}^{-1}E[p(x_{i})\gamma_{0}(x_{i})].
\]
Here we see that $E[\hat{\beta}_{H}-\check{\beta}|T]$ is the negative of a
Taylor expansion to order $Q$ of $\Lambda(\hat{\Sigma},\Sigma)$ in
$\hat{\Sigma}$ around $\Sigma.$ Therefore, it will follow that%
\[
E[\hat{\beta}_{H}-\beta_{0}|T]=O(K^{-(s_{\gamma}+s_{\alpha})/r})+O(\left\Vert
\hat{\Sigma}-\Sigma\right\Vert _{op}^{Q})=O(K^{-(s_{\gamma}+s_{\alpha}%
)/r})+O(\left(  \frac{K\ln(K)}{n}\right)  ^{Q/2}),
\]
where $\left\Vert \cdot\right\Vert _{op}$ is the operator norm for a matrix
and the second equality follows by the Rudelson (1999) matrix law of large
numbers. This equation is similar to the conclusion of Theorem 3 of Mukherjee,
Newey, and Robins (2017).

In comparison with the HOIF estimator the CF plug-in series estimator has a
remainder rate from estimating $\Sigma$ that is $\ln(n)K^{-s_{\gamma}/r}$ for
$s_{\gamma},s_{\alpha}\leq1$ and Haar splines and $\sqrt{K\ln(K)/n}%
K^{1/2-s_{\gamma}/r}$ more generally, without any higher order U-statistic
correction for the presence of $\hat{\Sigma}^{-1}.$ The DCDR estimator has
$\Delta_{n}=\Delta_{n}^{\ast},$ also without the need to rely on any
higher-order U-statistics. The key difference between the HOIF and these other
estimators is that the plug-in and doubly robust estimators use spline
regression in their construction and the HOIF estimator uses $\hat{\Sigma
}^{-1}$ from a training subsample.

Previously the HOIF estimator was the only known method of obtaining an
semiparametric efficient estimator of the expected conditional covariance
under the minimal conditions of Robins et al.(2009). We find here that the CF
plug-in estimator with a Haar basis can do this for $s_{\gamma},s_{\alpha}%
\leq1$ and for a general spline basis with $s_{\gamma}\geq r/2.$ We also find
that the DCDR estimator can do this for all $s_{\gamma}$ and $s_{\alpha}$.
These estimators are simpler than the HOIF\ estimator in not requiring the
higher order U-statistic terms$.$ It would be interesting to compare the size
of constants in respective remainder terms where HOIF could have an advantage
by virtue of its higher order influence function interpretation. That
comparison is beyond the scope of this paper.

The HOIF estimator remains the only known estimator that is semiparametric
efficient under the Robins et al.(2009) minimal conditions for the mean with
missing data over all $s_{\gamma}$ and $s_{\alpha}$. We expect that property
of HOIF to extend to all the linear average functionals we are considering in
this paper.

In summary, cross-fitting can be used to reduce bias of estimators and obtain
faster remainder rates. If cross fitting is not used for either the plug-in or
the doubly robust estimator there would be an additional $K/\sqrt{n}$ bias
term in the remainder. This extra term can increase the bias of the estimator
significantly for large $K.$ It is well known to be very important in some
settings, such as instrumental variables estimation as shown by Blomquist and
Dahlberg (1999) and Imbens, Angrist, and Krueger (1999). Also, its presence
prevents the plug-in estimator from attaining root-n consistency under minimal
conditions. Cross-fitting eliminates this large remainder for the linear
functionals we consider and results in plug-in and doubly robust estimators
with remainders that converge to zero as fast as known possible for
$s_{\gamma},s_{\alpha}\leq1,$ for $s_{\gamma}>r/2$, and for any $s_{\alpha}$
and $s_{\gamma}$ for a doubly robust estimator of the expected conditional covariance.

\section{Estimators of Average Linear Functionals}

We will analyze estimators of functionals of a conditional expectation%
\[
\gamma_{0}(x)=E[y_{i}|x_{i}=x],
\]
where $y_{i}$ is a scalar component and $x_{i}$ a subvector of $z_{i}$. Let
$\gamma$ represent a possible conditional expectation function and
$m(z,\gamma)$ denote a function of $\gamma$ and a possible realization $z$ of
a data observation. We consider
\[
\beta_{0}=E\left[  m(z_{i},\gamma_{0})\right]  ,
\]
where $m(z,\gamma)$ is an affine functional of $\gamma$ for every $z,$ meaning
$m(z,\gamma)-m(z,0)$ is linear in $\gamma$.

There are many important examples of such an object. One of these is the
expected conditional covariance we consider in Section 2. There $m(z,\gamma
)=a[y-\gamma(x)]$. This object shows up in different forms in the numerator
and denominator of%
\[
\xi_{0}=\frac{E[Cov(a_{i},y_{i}|x_{i})]}{E[Var(a_{i}|x_{i})]}.
\]
Here $\delta_{0}$ is the coefficient of $a_{i}$ in the population least
squares projection of $y_{i}$ on functions of the form $a_{i}\delta+g(x_{i}).$
Under an ignorability assumption this object $\delta_{0}$ can be interpreted
as a weighted average of conditional average treatment effects when $a_{i}$ is
a binary indicator for treatment and $x_{i}$ are covariates.

Another important example is the mean when data are missing at random. The
object of interest is $\beta_{0}=E[Y_{i}]$ where $Y_{i}$ is a latent variable
that is not always observed. Let $a_{i}$ be an observed binary indicator where
$a_{i}=1$ if $Y_{i}$ is observed. Suppose that there are observed covariates
$w_{i}$ such that $Y_{i}$ is mean independent of $a_{i}$ conditional on
$w_{i}$, i.e. $E[Y_{i}|a_{i}=1,w_{i}]=E[Y_{i}|w_{i}].$ Then for the observed
variable $y_{i}=a_{i}Y_{i}$ we have%
\[
E[E[y_{i}|a_{i}=1,w_{i}]]=E[E[Y_{i}|a_{i}=1,w_{i}]]=E[E[Y_{i}|w_{i}%
]]=\beta_{0}.
\]
Let $x=(a,w)$ and $\gamma_{0}(x_{i})=E[y_{i}|x_{i}].$ Then for $m(z,\gamma
)=\gamma(1,w)$ we have $\beta_{0}=E[m(z_{i},\gamma_{0})]$.

A third example is a weighted average derivative, where the object of interest
is%
\[
\beta_{0}=\int v(x)\left[  \partial\gamma_{0}(x)/\partial x_{1}\right]  dx,
\]
for some weight function $v(x),$ with $x_{1}$ continuously distributed and
$\int v(x)dx=1$. This object is proportional to $\beta_{10}$ in a conditional
mean index model where $E[y_{i}|x_{i}]=\tau(x_{i}^{T}\beta_{0})$ for some
unknown function $\tau(\cdot),$ as in Stoker (1986). This object is included
in the framework of this paper for $m(z,\gamma)=\int v(x)\left[
\partial\gamma(x)/\partial x_{1}\right]  dx.$ Assuming that $v(x)$ is zero at
the boundary, integration by parts gives%
\[
m(z,\gamma)=m(\gamma)=\int\omega(x)\gamma(x)dx,\omega(x)=-\partial
v(x)/\partial x_{1}.
\]

Throughout we will focus on the case where estimators of $\beta_{0}$ have a
finite semiparametric variance bound and so should be root-n consistently
estimable under sufficient regularity conditions. As discussed in Newey
(1994), this corresponds to $E[m\left(  z_{i},\gamma\right)  ]$ being mean
square continuous as a function of $\gamma$, so that by the Riesz
representation theorem the following condition is satisfied:

\bigskip

\textsc{Assumption 1:} \textit{There is }$\alpha_{0}\left(  x\right)
$\textit{ with }$E[\alpha_{0}(x_{i})^{2}]<\infty$\textit{ and for all }%
$\gamma$ with $E[\gamma(x_{i})^{2}]<\infty$,%
\begin{equation}
E\left[  m\left(  z_{i},\gamma\right)  -m(z_{i},0)\right]  =E\left[
\alpha_{0}\left(  x_{i}\right)  \gamma\left(  x_{i}\right)  \right]  .
\label{Riesz}%
\end{equation}

\bigskip

The function $\alpha_{0}(x)$ has an important role in the asymptotic theory.
The bias in a series estimator of $\beta_{0}$ will depend on the expected
product of biases in approximating $\gamma_{0}(x)$ and $\alpha_{0}(x)$.
Consequently there will be a trade-off in conditions that can be imposed on
$\gamma_{0}(x)$ and $\alpha_{0}(x)$ so that the estimators of $\beta_{0}$ have
good properties.

To help explain this condition we give the form of $\alpha_{0}(x)$ in each of
the examples. In the expected conditional covariance example iterated
expectations gives%
\begin{align}
E\left[  m\left(  z_{i},\gamma\right)  -m(z_{i},0)\right]   &  =-E[a_{i}%
\gamma(x_{i})]=-E[E[a_{i}|x_{i}]\gamma(x_{i})]=E[\alpha_{0}(x_{i})\gamma
(x_{i})],\label{ccriesz}\\
\alpha_{0}(x_{i})  &  =-E[a_{i}|x_{i}].\nonumber
\end{align}
In the missing data example, for the propensity score $\Pr(a_{i}=1|w_{i}%
)=\pi_{0}(w_{i})$, iterated expectations gives
\begin{align}
E\left[  m\left(  z_{i},\gamma\right)  -m(z_{i},0)\right]   &  =E[\gamma
(1,w_{i})]=E[\frac{\pi_{0}(w_{i})}{\pi_{0}(w_{i})}\gamma(1,w_{i}%
)]=E[\frac{a_{i}}{\pi_{0}(w_{i})}\gamma(1,w_{i})]\label{mdriesz}\\
&  =E[\frac{a_{i}}{\pi_{0}(w_{i})}\gamma(x_{i})]=E[\alpha_{0}(x_{i}%
)\gamma(x_{i})],\alpha_{0}(x_{i})=\frac{a_{i}}{\pi_{0}(w_{i})}.\nonumber
\end{align}
In the average derivative example, multiplying and dividing by the pdf
$f_{0}(x)$ of $x_{i}$ gives%
\begin{align}
E\left[  m\left(  z_{i},\gamma\right)  -m(z_{i},0)\right]   &  =\int%
\omega(x)\gamma(x)dx=\int\frac{\omega(x)}{f_{0}(x)}\gamma(x)f_{0}%
(x)dx=E[\frac{\omega(x_{i})}{f_{0}(x_{i})}\gamma(x_{i})]\label{adriesz}\\
&  =E[\alpha_{0}(x_{i})\gamma(x_{i})],\alpha_{0}(x_{i})=\frac{\omega(x_{i}%
)}{f_{0}(x_{i})}.\nonumber
\end{align}

Our estimators of $\beta_{0}$ will be based on a nonparametric estimator
$\hat{\gamma}$ of $\gamma_{0}$ and possibly on a nonparametric estimator
$\tilde{\alpha}$ of $\alpha_{0}.$ The CF plug-in estimator is given by%
\[
\hat{\beta}=\frac{1}{n}\sum_{\ell=1}^{L}\sum_{i\in I_{\ell}}m(z_{i}%
,\hat{\gamma}_{\ell}),
\]
where $I_{\ell},(\ell=1,...,L)$ is a partition of the observation index set
$\{1,...,n\}$ into $L$ distinct subsets of about equal size and $\hat{\gamma
}_{\ell}$ only uses observations \textit{not} in $I_{\ell}.$ We will consider
a fixed number of groups $L$ in the asymptotics. It would be interesting to
consider results where the number of groups grows with the sample size, even
"leave one out" estimators where $I_{\ell}$ only includes one observation, but
theory for those estimators is more challenging and we leave it to future work.

The DCDR estimator makes use of $\tilde{\alpha}_{\ell}$ that may be
constructed from different observations than $\hat{\gamma}_{\ell}.$ The doubly
robust estimator is%
\[
\tilde{\beta}=\frac{1}{n}\sum_{\ell=1}^{L}\sum_{i\in I_{\ell}}\left\{
m(z_{i},\hat{\gamma}_{\ell})+\tilde{\alpha}_{\ell}(x_{i})[y_{i}-\hat{\gamma
}_{\ell}(x_{i})]\right\}  .
\]
This estimator has the form of a plug-in estimator plus the sample average of
$\tilde{\alpha}_{\ell}(x_{i})[y_{i}-\hat{\gamma}_{\ell}(x_{i})],$ which is an
estimator of the influence function of $\int m(z,\hat{\gamma}_{\ell}%
)F_{0}(dz).$ The addition of $\tilde{\alpha}_{\ell}(x_{i})[y_{i}-\hat{\gamma
}_{\ell}(x_{i})]$ will mean that the nonparametric estimators $\hat{\gamma
}_{\ell}$ and $\tilde{\alpha}_{\ell}$ do not affect the asymptotic
distribution of $\tilde{\beta},$ i.e. the limit distribution would be the same
if $\hat{\gamma}_{\ell}$ and $\tilde{\alpha}_{\ell}$ were replaced by their
true values and $\Delta_{n}\longrightarrow0$. This estimator allows for full
cross-fitting where $\tilde{\alpha}$ and $\hat{\gamma}$ may be based on
distinct subsamples.

The cross-fit estimator $\tilde{\beta}$ is doubly robust in the sense that
$\tilde{\beta}$ will be consistent as long as either $\hat{\gamma}_{\ell}$ or
$\tilde{\alpha}_{\ell}$ is consistent, as shown by Chernozhukov et al.(2016)
for this general class of functionals. When $\hat{\gamma}(x)$ is a series
estimator like that described above the CF plug-in estimator $\hat{\beta}$ is
also doubly robust in a more limited sense. It will be consistent with fixed
$p(x)$ if either $\gamma_{0}(x)$ or $\alpha_{0}(x)$ is a linear combination of
$p(x)$, as shown for the mean with missing data in Robins et al.(2007) and in
Chernozhukov et al.(2016) for the general linear function case we are considering.

Throughout the paper we assume that each data point $z_{i}$ is used for
estimation for some group $\ell$ and that the number of observations in group
$\ell$, the number used to form $\hat{\gamma}_{\ell}$, and the number used to
form $\tilde{\alpha}_{\ell}$ grow at the same rate as the sample size. To make
this condition precise let $\bar{n}_{\ell}$ be the number of elements in
$I_{\ell},$ $\hat{n}_{\ell}$ be the number used to form $\hat{\gamma}_{\ell},$
and $\tilde{n}_{\ell}$ be the number of observations used to form
$\tilde{\alpha}_{\ell}$. We will assume throughout that all the observations
are used for each $\ell$, i.e. that either $\bar{n}_{\ell}+\hat{n}_{\ell}=n$
or $\bar{n}_{\ell}+\hat{n}_{\ell}+\tilde{n}_{\ell}=n$ if different
observations are used for $\hat{\gamma}_{\ell}$ and $\tilde{\alpha}_{\ell}$.

\bigskip

\textsc{Assumption 2:} \textit{ There is a constant }$C>0$\textit{ such that
either} $\bar{n}_{\ell}+\hat{n}_{\ell}=n$\textit{ and }$\min_{\ell}\{\bar
{n}_{\ell},\hat{n}_{\ell}\}\geq Cn$\textit{ or }$\bar{n}_{\ell}+\hat{n}_{\ell
}+\tilde{n}_{\ell}=n$\textit{ and }$\min_{\ell}\{\bar{n}_{\ell},\hat{n}_{\ell
},\tilde{n}_{\ell}\}\geq Cn.$ \textit{For the plug-in estimator groups are as
close as possible to being of equal size.}

\bigskip

The assumption that the group sizes are as close to equal as possible for the
plug-in estimator is made for simplicity but could be relaxed.

We turn now to conditions for the regression spline estimators of $\gamma
_{0}(x)$ and $\alpha_{0}(x)$. We continue to consider regression spline first
steps where $p(x)$ is a $K\times1$ vector of regression splines. The
nonparametric estimator of $\gamma_{0}(x)$ will be a series regression
estimator where%
\[
\hat{\gamma}_{\ell}(x)=p(x)^{T}\hat{\delta}_{\ell},\text{ }\hat{\delta}_{\ell
}=\hat{\Sigma}_{\ell}^{-}\hat{h}_{\ell},\text{ }\hat{\Sigma}_{\ell}=\frac
{1}{\hat{n}_{\ell}}\sum_{i\in\hat{I}_{\ell}}p(x_{i})p(x_{i})^{T},\text{ }%
\hat{h}_{\ell}=\frac{1}{\hat{n}_{\ell}}\sum_{i\in\hat{I}_{\ell}}p(x_{i}%
)y_{i},
\]
where a $T$ superscript denotes the transpose, $\hat{I}_{\ell}$ is the index
set for observations used to construct $\hat{\gamma}_{\ell}(x)$, and $A^{-}$
denotes any generalized inverse of a positive semi-definite matrix $A$. Under
conditions given below $\hat{\Sigma}_{\ell}$ will be nonsingular with
probability approaching one so that $\hat{\Sigma}_{\ell}^{-}=\hat{\Sigma
}_{\ell}^{-1}$ for each $\ell.$

The DCDR estimator $\tilde{\beta}$ uses an estimator of $\alpha_{0}(x).$ The
function $\alpha_{0}(x)$ cannot generally be interpreted as a conditional
expectation and so cannot generally be estimated by a linear regression.
Instead we use Assumption 1 and equation (\ref{Riesz}) to construct an
estimator. Let $v(z)=(m(z,p_{1})-m(z,0),...,m(z,p_{K})-m(z,0))^{T}$. Then by
Assumption 1,%
\[
E[v(z_{i})]=E[p(x_{i})\alpha_{0}(x_{i})],
\]
so that $\tilde{h}_{\ell\alpha}=\sum_{i\in\tilde{I}_{\ell}}v(z_{i})/\tilde
{n}_{\ell}$ is an unbiased estimator of $E[p(x_{i})\alpha_{0}(x_{i})].$ A
series estimator of $\alpha_{0}(x)$ is then%
\[
\tilde{\alpha}_{\ell}(x)=p(x)^{T}\tilde{\delta}_{\ell\alpha},\tilde{\delta
}_{\ell\alpha}=\tilde{\Sigma}_{\ell}^{-}\tilde{h}_{\ell\alpha},\tilde{\Sigma
}_{\ell}=\frac{1}{\tilde{n}_{\ell}}\sum_{i\in\tilde{I}_{\ell}}p(x_{i}%
)p(x_{i})^{T}.
\]
Here $\tilde{\delta}_{\ell\alpha}$ is an estimator of the coefficients of the
population regression of $\alpha_{0}(x)$ on $p(x),$ but $\tilde{\delta}%
_{\ell\alpha}$ is not obtained from a linear regression. This type of
estimator of $\alpha_{0}(x)$ was used to construct standard errors for
functionals of series estimators in Newey (1994).

Now that we have specified the form of the estimators $\hat{\gamma}_{\ell}$
and $\tilde{\alpha}_{\ell}$ we can give a complete description of the
estimators in each of the examples. For the expected conditional covariance
recall that $m(z,\gamma)=a[y-\gamma(x)].$ Therefore the CF plug-in estimator
will be%
\begin{equation}
\hat{\beta}=\frac{1}{n}\sum_{\ell=1}^{L}\sum_{i\in I_{\ell}}a_{i}[y_{i}%
-\hat{\gamma}_{\ell}(x_{i})]. \label{ccplug}%
\end{equation}
Also, as discussed above, for the expected conditional covariance $\alpha
_{0}(x)=-E[a_{i}|x_{i}=x]$ and $v(z_{i})=-a_{i}p(x_{i})$, so that
$\tilde{\alpha}_{\ell}(x)=-\tilde{\gamma}_{a\ell}(x)$ where $\tilde{\gamma
}_{a\ell}(x)=p(x)^{T}\tilde{\Sigma}_{\ell}^{-}\sum_{i\in\tilde{I}_{\ell}%
}p(x_{i})a_{i}/\tilde{n}_{\ell}$ is the regression of $a_{i}$ on $p(x_{i})$
for the observations indexed by $\tilde{I}_{\ell}.$ Then the DCDR estimator is%
\begin{equation}
\tilde{\beta}=\frac{1}{n}\sum_{\ell=1}^{L}\sum_{i\in I_{\ell}}[a_{i}%
+\tilde{\alpha}_{\ell}(x_{i})][y_{i}-\hat{\gamma}_{\ell}(x_{i})]=\frac{1}%
{n}\sum_{\ell=1}^{L}\sum_{i\in I_{\ell}}\{a_{i}-\tilde{E}[a_{i}|x_{i}%
]\}[y_{i}-\hat{\gamma}_{\ell}(x_{i})], \label{ccrob}%
\end{equation}
where $\tilde{E}[a_{i}|x_{i}]=-\tilde{\alpha}_{\ell}(x_{i})$ is the predicted
value from the regression of $a_{i}$ on $p(x_{i}).$ This estimator is the
average of the product of two nonparametric regression residuals, where the
average and each of the nonparametric estimators can be constructed from
different samples.

For the missing data example the estimators are based on series estimation of
$E[y_{i}|a_{i}=1,w_{i}]$. Let $q(w)$ denote a $K\times1$ vector of splines,
$x=(a,w^{T})^{T},$ and $p(x)=(aq(w)^{T},(1-a)q(w)^{T})^{T}$. The predicted
value $\hat{\gamma}(1,w)$ will be the same as from a linear regression of
$y_{i}$ on $q(w_{i})$ for observations with $a_{i}=1.$ That is, $\hat{\gamma
}(1,w)=q(w)^{T}\hat{\delta}_{\ell}$ where
\[
\hat{\delta}_{\ell}=\hat{\Sigma}_{\ell}^{-}\hat{h}_{\ell},\text{ }\hat{\Sigma
}_{\ell}=\frac{1}{\hat{n}_{\ell}}\sum_{i\in\hat{I}_{\ell}}a_{i}q(w_{i}%
)q(w_{i})^{T},\text{ }\hat{h}_{\ell}=\frac{1}{\hat{n}_{\ell}}\sum_{i\in\hat
{I}_{\ell}}a_{i}q(w_{i})y_{i}.
\]
The CF plug-in estimator is%
\[
\hat{\beta}=\frac{1}{n}\sum_{\ell=1}^{L}\sum_{i\in I_{\ell}}q(w_{i})^{T}%
\hat{\delta}_{\ell}.
\]
The DCDR estimator is based on an estimator of the inverse propensity score
$\pi_{0}(w_{i})^{-1}=1/\pi_{0}(w_{i})$ given by%
\[
\widetilde{\pi(w_{i})_{\ell}^{-1}}=q(w_{i})^{T}\tilde{\delta}_{\ell}^{\alpha
},\tilde{\delta}_{\ell}^{\alpha}=\tilde{\Sigma}_{\ell}^{-}\tilde{h}_{\ell
}^{\alpha},\tilde{\Sigma}_{\ell}=\frac{1}{\tilde{n}_{\ell}}\sum_{i\in\tilde
{I}_{\ell}}a_{i}q(w_{i})q(w_{i})^{T},\text{ }\tilde{h}_{\ell}^{\alpha}%
=\frac{1}{\tilde{n}_{\ell}}\sum_{i\in\tilde{I}_{\ell}}q(w_{i}),
\]
where $\tilde{n}_{\ell}$ is the number of observation indices in $\tilde
{I}_{\ell}$. This estimator of the inverse propensity score is a version of
one discussed in Robins et al.(2007). The DCDR estimator is%
\[
\tilde{\beta}=\frac{1}{n}\sum_{\ell=1}^{L}\sum_{i\in I_{\ell}}\left\{
q(w_{i})^{T}\hat{\delta}_{\ell}+a_{i}\widetilde{\pi(w_{i})_{\ell}^{-1}}%
[y_{i}-q(w_{i})^{T}\hat{\delta}_{\ell}]\right\}  .
\]
This has the usual form for a doubly robust estimator of the mean with data
missing at random. It differs from previous estimators in having the full CF
form where the nonparametric estimators are based on distinct subsamples of
the data.

For the average derivative example $m(z,\gamma)=\int\omega(x)\gamma(x)dx$ does
not depend on $z$ so we can use all the data in the construction of the
plug-in estimator. That estimator is given by%
\begin{equation}
\hat{\beta}=\int\omega(x)\hat{\gamma}(x)dx=v^{T}\hat{\delta}\text{, }%
v=\int\omega(x)p(x)dx,\hat{\delta}=[\sum_{i=1}^{n}p(x_{i})p(x_{i})^{T}%
]^{-}\sum_{i=1}^{n}p(x_{i})y_{i}. \label{adplug}%
\end{equation}
As shown in equation (\ref{adriesz}), $\alpha_{0}(x)=f_{0}(x)^{-1}\omega(x),$
where $f_{0}(x)$ is the pdf of $x$. Also here $v(z)=v$ so the estimator of
$\alpha_{0}(x)$ is $p(x)^{T}\tilde{\Sigma}_{\ell}^{-}v.$ The DCDR estimator is
then
\begin{equation}
\tilde{\beta}=\frac{1}{n}\sum_{\ell=1}^{L}\sum_{i\in I_{\ell}}\{\int%
\omega(x)\hat{\gamma}_{\ell}(x)dx+\left[  p(x_{i})^{T}\tilde{\Sigma}_{\ell
}^{-}v\right]  [y_{i}-\hat{\gamma}_{\ell}(x_{i})]\}. \label{adrob}%
\end{equation}
Both the plug-in and the DCDR estimators depend on the integral $v=\int%
\omega(x)p(x)dx.$ Generally this vector of integrals will not exist in closed
form so that construction of these estimators will require numerical
computation or estimation of $v$, such as by simulation.

We now impose some specific conditions on $p(x)$.

\bigskip

\textsc{Assumption 3:} $p(x)=aq(w)$\textit{ where i) the support of }$w_{i}%
$\textit{ is }$[0,1]^{r}$\textit{, }$w_{i}$\textit{ is continuously
distributed with bounded pdf that is bounded away from zero; ii) }%
$q(w)$\textit{ are tensor product b-splines of order }$\kappa$\textit{ with
knot spacing approximately proportional to the number of knots; iii) }%
$q(w)$\textit{ is normalized so that }$\lambda_{\min}(E[q(w_{i})q(w_{i}%
)^{T}])\geq C>0$\textit{ and }$\sup_{w\in\lbrack0,1]^{r}}\left\Vert
q(w)\right\Vert \leq C\sqrt{K};$\textit{ iv) }$a_{i}$\textit{ is bounded and
}$E[a_{i}^{2}|w_{i}]$\textit{ is bounded away from zero.}

\bigskip

Under condition i) it is known that there is a normalization such that
condition iii) is satisfied, e.g. as in Newey (1997). To control the bias of
the estimator we require that the true regression function $\gamma_{0}(x)$ and
the auxiliary function $\alpha_{0}(x)$ each be in a Holder class of functions.
We define a function $g(x)$ to be Holder of order $s$ if there is a constant
$C$ such that $g(x)$ is continuously differentiable of order $\bar{s}=int[s]$
and each of its $\bar{s}$ partial derivatives $\nabla^{\bar{s}}g(x)$ satisfies
$\left\vert \nabla^{\bar{s}}g(\tilde{x})-\nabla^{\bar{s}}g(x)\right\vert \leq
C\left\Vert \tilde{x}-x\right\Vert ^{s-\bar{s}}.$

$\bigskip$

\textsc{Assumption 4:} $\gamma_{0}(x)$\textit{ and }$\alpha_{0}(x)$\textit{
are Holder of order }$s_{\gamma}$\textit{ and }$s_{\alpha}$\textit{
respectively.}

\bigskip

This condition implies that the population least squares approximations to
$\gamma_{0}(x)$ and \textit{ }$\alpha_{0}(x)$ converge at certain rates. Let
$\zeta_{\gamma}=\min\{1+\kappa,s_{\gamma}\}/r,$ $\zeta_{\alpha}=\min
\{1+\kappa,s_{a}\}/r,$ $\Sigma=E[p(x_{i})p(x_{i})^{T}]$, $\delta=\Sigma
^{-1}E[p(x_{i})\gamma_{0}(x_{i})],$ $\gamma_{K}(x)=p(x)^{T}\delta,$
$\delta_{\alpha}=\Sigma^{-1}E[p(x_{i})\alpha_{0}(x_{i})],$ $\alpha
_{K}(x)=p(x)^{T}\delta_{a}.$ Then standard approximation theory for splines
gives%
\begin{align*}
E[\{\gamma_{0}(x_{i})-\gamma_{K}(x_{i})\}^{2}]  &  =O(K^{-2\zeta_{\gamma}%
}),\sup_{x\in\lbrack0,1]^{r}}\left\vert \gamma_{0}(x)-\gamma_{K}(x)\right\vert
=O(K^{-\zeta_{\gamma}}),\\
E[\{\alpha_{0}(x_{i})-\alpha_{K}(x_{i})\}^{2}]  &  =O(K^{-2\zeta_{\alpha}}).
\end{align*}
We will use these results to derive the rates at which certain remainders
converge to zero.

We also impose the following condition:

\bigskip

\textsc{Assumption 5: }$Var(y_{i}|x_{i})\leq C,$ $K\longrightarrow\infty$,
\textit{and} $K\ln(K)/n\longrightarrow0$.

\bigskip

These are standard conditions for series estimators of conditional
expectations. A bounded conditional variance for $y_{i}$ helps bound the
variance of series estimators. The upper bound on the rate at which $K$ grows
is slightly stronger than $K/n\longrightarrow0$. This upper bound on $K$
allows us to apply the Rudelson (1999) law of large numbers for symmetric
matrices to show that the various second moment matrices of $p(x)$ converge in
probability. Another condition we impose is:

\bigskip

\textsc{Assumption 6:} $\lambda_{\max}(E[v(z_{i})v(z_{i})^{T}])\leq Cd_{K}$
\textit{and }$\left\{  E[\{m(z_{i},\gamma_{K})-m(z_{i},\gamma_{0}%
)\}^{2}]\right\}  ^{1/2}=O(K^{-\zeta_{m}}).$

\bigskip

The first condition will be satisfied with $d_{K}=1$ in the examples under
specific regularity conditions detailed below. The second condition gives a
rate for the mean square error convergence of $m(z,\gamma_{K})-m(z,\gamma
_{0})$ as $K$ grows. In all of the examples this rate will be $\zeta_{m}%
=\zeta_{\gamma}.$ In other examples, including those where $m(z,\gamma)$ and
$v(z)$ depend on derivatives with respect to $x,$ we will have $d_{K}$ growing
with $K$ and $\zeta_{m}<\zeta_{\gamma}.$

For the statement of the results to follow it is convenient to work with the
remainder term%
\[
\bar{\Delta}_{n}^{\ast}=\sqrt{n}K^{-\zeta_{\gamma}-\zeta_{\alpha}}%
+K^{-\zeta_{\gamma}}+K^{-\zeta_{\alpha}}+\sqrt{\frac{K}{n}}.
\]
This remainder coincides with the fast remainder $\Delta_{n}^{\ast}$ when the
spline order is high enough with $\kappa\geq\max\{s_{\gamma},s_{\alpha}\}-1.$
The only cases where it would not be possible to choose such a $\kappa$ are
for the Haar basis where $\kappa=0.$

\bigskip

\section{The Plug-in Estimator}

In this Section we derive bounds on the size of remainders for the plug-in
estimator. Some bounds are given for general plug-in estimators, some for
plug-ins that are series regression with Haar splines, and some for other
splines. We begin with a result that applies to all plug-ins. We drop the CF
designation because all the estimators from this point on will use cross-fitting.

The cross-fit form of the plug-in estimator allows us to partly characterize
its properties under weak conditions on a general plug-in estimator that need
not be a series regression. This characterization relies on independence of
$\hat{\gamma}_{\ell}$ from the observations in $I_{\ell}$ to obtain relatively
simple stochastic equicontinuity remainders. Also, this result accounts for
the overlap across groups in observations used to form $\hat{\gamma}_{\ell}$.
Let $\mathcal{A}_{n}$ denote an event that occurs with probability approaching
one. For example, $\mathcal{A}_{n}$ could include the set of data points where
$\hat{\Sigma}_{\ell}$ is nonsingular for each $\ell.$

\bigskip

\textsc{Lemma 1:} \textit{If Assumptions 1 and 2 are satisfied and there is
}$\Delta_{n}^{m}$ \textit{such that }%
\[
1(\mathcal{A}_{n})\left\{  \int[m(z,\hat{\gamma}_{\ell})-m(z,\gamma_{0}%
)]^{2}F_{0}(dz)\right\}  ^{1/2}=O_{p}(\Delta_{n}^{m}),\left(  \ell
=1,...,L\right)  ,
\]
\textit{then for }$\bar{m}(\gamma)=\int m(z,\gamma)F_{0}(dz),$%
\[
\sqrt{n}(\hat{\beta}-\beta_{0})=\frac{1}{\sqrt{n}}\sum_{i=1}^{n}%
[m(z_{i},\gamma_{0})-\beta_{0}]+\sqrt{n}\sum_{\ell=1}^{L}\frac{\bar{n}_{\ell}%
}{n}[\bar{m}(\hat{\gamma}_{\ell})-\beta_{0}]+O_{p}(\Delta_{n}^{m}).
\]
\textit{If in addition there is }$\Delta_{n}^{\phi}$ \textit{such that for
each }$\left(  \ell=1,...,L\right)  ,$
\[
\sqrt{\hat{n}_{\ell}}[\bar{m}(\hat{\gamma}_{\ell})-\beta_{0}]=\frac{1}%
{\sqrt{\hat{n}_{\ell}}}\sum_{i\notin I_{\ell}}\alpha_{0}(x_{i})[y_{i}%
-\gamma_{0}(x_{i})]+O_{p}(\Delta_{n}^{\phi}),
\]
\textit{then for }$\delta(z)=m(z,\beta_{0})-\beta_{0}+\alpha_{0}%
(x)[y-\gamma_{0}(x)]$%
\[
\sqrt{n}(\hat{\beta}-\beta_{0})=\frac{1}{\sqrt{n}}\sum_{i=1}^{n}\psi
(z_{i})+O_{p}(\Delta_{n}^{m}+\Delta_{n}^{\phi}+n^{-1}).
\]
\textit{ }

\bigskip

This result gives a decomposition of remainder bounds into two kinds. The
first $\Delta_{n}^{m}$ is a stochastic equicontinuity bound that has the
simple mean-square form given here because of the sample splitting. The second
$\Delta_{n}^{\phi}$ is a bound that comes from the asymptotically linear
expansion of the linear functional estimator $\bar{m}(\hat{\gamma}_{\ell})$.
For general b-splines we can apply Ichimura and Newey (2017) to obtain
$\Delta_{n}^{\phi}$. For zero order splines we give here sharper remainder bounds.

For series estimators the stochastic equicontinuity remainder bound
$\Delta_{n}^{m}$ will be%
\[
\Delta_{n}^{m}=\sqrt{(d_{K}+1)\frac{K}{n}}+K^{-\zeta_{m}},
\]
where $d_{K}$ and $\zeta_{m}$ are as given in Assumption 6. As mentioned
above, in the examples in this paper $d_{K}\leq C$ and $\zeta_{m}%
=\zeta_{\gamma}$. Here we can take $\Delta_{n}^{m}\leq C\bar{\Delta}_{n}%
^{\ast}$, so the stochastic equicontinuity remainder bound is the same size as
$\bar{\Delta}_{n}^{\ast}$.

Our next result gives remainder bounds for the Haar basis.

\bigskip

\textsc{Theorem 2:} \textit{If Assumptions 1-6 are satisfied, }$\kappa
=0,$\textit{ and }$K[\ln(n)]^{2}/n\longrightarrow0$ \textit{then}%
\[
\sqrt{n}(\hat{\beta}-\beta_{0})=\frac{1}{\sqrt{n}}\sum_{i=1}^{n}\psi
(z_{i})+O_{p}(\bar{\Delta}_{n}^{\ast}+\Delta_{n}^{m}+K^{-\zeta_{\gamma}}%
\ln(n)).
\]
\textit{If in addition }$d_{K}$\textit{ is bounded as a function of }%
$K$\textit{ and }$\zeta_{m}=\zeta_{\gamma}$\textit{ then }$\Delta_{n}^{m}\leq
C\bar{\Delta}_{n}^{\ast}$\textit{.}

\bigskip

Here we see that for a Haar basis the order of the remainder term for the
plug-in estimator is a sum of the stochastic equicontinuity term $\Delta
_{n}^{m}$ and $\bar{\Delta}_{n}^{\ast}$, with $K^{-\zeta_{\gamma}}\ln(n)$
being the size of the fast remainder up to $\ln(n).$ In the examples and other
settings where $d_{K}$ is bounded and $\zeta_{m}=\zeta_{\gamma}$ the
$\Delta_{n}^{m}$ remainder will just be of order $\bar{\Delta}_{n}^{\ast}$.
The following result states conditions for the examples.

\bigskip

\textsc{Corollary 3:} \textit{Suppose that Assumptions 1-3 and 5 are
satisfied, }$\kappa=0$\textit{, }$K[\ln(n)]^{2}/n\longrightarrow0,$
\textit{and }$\gamma_{0}(x)$\textit{ is Holder of order }$s_{\gamma}.$\textit{
If either i) }$\hat{\beta}$\textit{ is the expected conditional covariance
estimator, }$E[a_{i}|x_{i}=x]$\textit{ is Holder of order }$s_{\alpha}%
$\textit{, }$E[a_{i}^{2}|x_{i}]$\textit{ is bounded, or ii) }$\hat{\beta}%
$\textit{ is the missing data mean estimator, }$\Pr(a_{i}=1|x_{i})$\textit{ is
bounded away from zero and is Holder of order }$s_{\alpha},$\textit{ or iii)
}$\hat{\beta}$\textit{ is the average derivative estimator, }$\omega
(x)$\textit{ and }$f_{0}(x)$\textit{ are Holder of order }$s_{\alpha}%
$\textit{, and }$f_{0}(x)$\textit{ is bounded away from zero on the set where
}$\omega(x)>0,$\textit{ then}%
\[
\sqrt{n}(\hat{\beta}-\beta_{0})=\frac{1}{\sqrt{n}}\sum_{i=1}^{n}\psi
(z_{i})+O_{p}(\bar{\Delta}_{n}^{\ast}+K^{-\zeta_{\gamma}}\ln(n)).
\]

\bigskip

The remainder bound means that the plug-in estimator can attain root-n
consistency under minimal conditions, when the dimension $r$ is small enough.
There will exist $K$ such that $\bar{\Delta}_{n}^{\ast}$ goes to zero if an
only if
\begin{equation}
1/2<\zeta_{\gamma}+\zeta_{\alpha}=\frac{\min\{1,s_{\gamma}\}+\min
\{1,s_{\alpha}\}}{r}. \label{min cond}%
\end{equation}
This condition can be satisfied for $r<4$ but not for $r\geq4.$ For $r=1$ this
condition will be satisfied if and only if
\[
s_{\gamma}+s_{\alpha}>\frac{1}{2},
\]
which is the minimal condition of Robins et al.(2009) for existence of a
semiparametric efficient estimator for the expected conditional covariance and
missing data parameters when $r=1$. For $r=2$ we note that%
\[
\min\{1,s_{\gamma}\}+\min\{1,s_{\alpha}\}\geq1\text{ if and only if }%
s_{\gamma}+s_{\alpha}\geq1.
\]
For $r=2$ equation (\ref{min cond}) is $\min\{1,s_{\gamma}\}+\min
\{1,s_{\alpha}\}>1$, which requires both $s_{\alpha}>0$ and $s_{\gamma}>0$ and
so is slightly stronger than the Robins et al.(2009) condition $s_{\gamma
}+s_{\alpha}>1$. For $r=3$ the situation is more complicated. Equation
(\ref{min cond}) is stronger than the corresponding condition $s_{\gamma
}+s_{\alpha}>3/2$ of Robins et al.(2009), although it is the same for the set
of $(s_{\gamma},s_{\alpha})$ where $s_{\gamma}\leq1$ and $s_{\alpha}\leq1.$
Along the diagonal where $s_{\alpha}=s_{\gamma}$ the two conditions coincide
as $s_{\gamma}>3/4.$

The limited nature of these results is associated with the Haar basis, which
limits the degree to which smoothness of the underlying function results in a
faster approximation rate. If Theorem 2 and Corollary 3 could be extended to
other, higher order b-splines, this limitation could be avoided. For the
present we are only able to do this for the doubly robust estimator of a
partially linear projection, as discussed in the next Section.

There is a key result that allows us to obtain the remainder bound
$\bar{\Delta}_{n}^{\ast}$ in Theorem 2. Let $\hat{h}_{2}=\sum_{i=1}^{n}%
p(x_{i})[\gamma_{0}(x_{i})-\gamma_{K}(x_{i})]/n$, $\hat{\Sigma}=\sum_{i=1}%
^{n}p(x_{i})p(x_{i})^{T}/n$, and $\Sigma=E[p(x_{i})p(x_{i})^{T}.$ We show in
the Appendix that for the Haar basis
\begin{equation}
\lambda_{\max}(E[(\Sigma-\hat{\Sigma})^{j}\hat{h}_{2}\hat{h}_{2}^{T}%
(\Sigma-\hat{\Sigma})^{j}])\leq\frac{K^{-2\zeta_{\gamma}}}{n}\left(  \frac
{CK}{n}\right)  ^{j}. \label{key result}%
\end{equation}
If b-spline bases other than Haar also satisfied this condition then we could
obtain results analogous to Theorem 2 and Corollary 3 for these bases. We do
not yet know if other bases satisfy this condition. The Haar basis is
convenient in $p(x)^{T}p(x)$ being piecewise constant. Cattaneo and Farrell
(2013) exploited other special properties of the Haar basis to obtain sharp
uniform nonparametric rates.

For b-splines of any order we can obtain remainder rates by combining Lemma 1
with Theorem 8 of Ichimura and Newey (2017).

\bigskip

\textsc{Theorem 4:} \textit{If Assumptions 1-6 are satisfied} \textit{then}%
\[
\sqrt{n}(\hat{\beta}-\beta_{0})=\frac{1}{\sqrt{n}}\sum_{i=1}^{n}\psi
(z_{i})+O_{p}(\bar{\Delta}_{n}^{\ast}+\Delta_{n}^{m}+\bar{\Delta}_{n}%
),\bar{\Delta}_{n}=\left(  \frac{K\ln K}{n}\right)  ^{1/2}K^{(1/2)-\zeta
_{\gamma}}.
\]
\textit{If in addition }$d_{K}$\textit{ is bounded as a function of }%
$K$\textit{ and }$\zeta_{m}=\zeta_{\gamma}$\textit{ then }$\Delta_{n}^{m}\leq
C\bar{\Delta}_{n}^{\ast}$\textit{.}

\bigskip

Here we see that the remainder bound for splines with $\kappa>0$ has an
additional term $\bar{\Delta}_{n}$. When $\zeta_{\gamma}$ is large enough,
i.e. $\gamma_{0}(x)$ is smooth enough and the order of the spline is big
enough, so that $\zeta_{\gamma}>1/2,$ the additional $\bar{\Delta}_{n}$ will
be no larger than $\bar{\Delta}_{n}^{\ast}.$ Also, when $\zeta_{\gamma}>1/2$
the condition of Robins et al.(2009) for semiparametric efficient estimation
is met for the expected conditional covariance and missing data examples for
any $\zeta_{\alpha}$. Thus, when $\gamma_{0}(x)$ is smooth enough to meet the
Robins et al.(2009) condition without imposing any smoothness on $\alpha
_{0}(x)$ the plug-in estimator will have the remainder bound $\bar{\Delta}%
_{n}^{\ast}.$

More generally there will exist a $K$ such that $\bar{\Delta}_{n}+\bar{\Delta
}_{n}^{\ast}$ goes to zero if and only if
\begin{equation}
2\min\{\kappa+1,s_{\gamma}\}+\min\{\kappa+1,s_{\alpha}\}>r. \label{minimal}%
\end{equation}
This condition is slightly stronger than that of Robins et al.(2009) which is
$2s_{\gamma}+2s_{\alpha}>r.$ Also, the remainder may go to zero when when $K$
is chosen to maximize the rate at which the mean square error of $\hat{\gamma
}_{0}(x)$ goes to zero. Setting $K^{-2\zeta_{r}}$ proportional to $K/n$ is
such a choice of $K$. Here the remainder term goes to zero for $\min
\{\kappa+1,s_{\gamma}\}>r/\left[  2(1+r)\right]  $ and $\min\{\kappa
+1,s_{\alpha}\}>r/2,$ a stronger condition for $s_{\gamma}$ and the same
condition for $s_{\alpha}$ as would hold if the remainder were $\bar{\Delta
}_{n}^{\ast}$.

\section{Partially Linear Projection}

In this Section we consider a series estimator of partially linear projection
coefficients. We give this example special attention because the DCDR
estimator will have a remainder bound that is only $\bar{\Delta}_{n}^{\ast}$.
The remainder bounds we find for other doubly robust estimators may be larger.
What appears to make the partially linear projection special in this respect
is that $\alpha_{0}(x)$ is a conditional expectation of an observed variable.
In other cases where $\alpha_{0}(x)$ is not a conditional expectation we do
not know if the remainder bound will be $\bar{\Delta}_{n}^{\ast}$ for bases
other than Haar.

The parameter vector of interest in this Section is%
\[
\beta_{0}=\left(  E[\{a_{i}-E[a_{i}|w_{i}]\}a_{i}^{T}]\right)  ^{-1}%
E[\{a_{i}-E[a_{i}|w_{i}]\}y_{i}].
\]
This vector $\beta_{0}$ can be thought of as the coefficients of $a_{i}$ in a
projection of $y_{i}$ on the set of functions of the form $a_{i}^{T}%
\beta+\lambda(x_{i})$ that have finite mean square. Note that this definition
of $\beta_{0}$ places no substantive restrictions on the distribution of data,
unlike the conditional expectation partially linear model where $E[y_{i}%
|a_{i},w_{i}]=a_{i}^{T}\beta_{0}+\xi_{0}(x_{i}).$

The object $\beta_{0}$ is of interest in a treatment effects model where
$a_{i}$ is a binary treatment, $y_{i}$ is the observed response, $x_{i}$ are
covariates, and outcomes with and without treatment are assumed to be mean
independent of $a_{i}$ conditional on $w_{i}$. Under an ignorability condition
that the outcome is mean independent of treatment conditional on covariates,
$E[y_{i}|a_{i}=1,x_{i}]-E[y_{i}|a_{i}=0,x_{i}]$ is the average treatment
effect conditional on $x_{i}$. Also for $\pi_{i}=\Pr(a_{i}=1|x_{i}),$%
\[
\beta_{0}=\frac{E[\pi_{i}(1-\pi_{i})\{E[y_{i}|a_{i}=1,x_{i}]-E[y_{i}%
|a_{i}=0,x_{i}]\}]}{E[\pi_{i}(1-\pi_{i})]}.
\]
Here we have the known interpretation of $\beta_{0}$ as a weighted average of
conditional average treatment effects, with weights $\pi_{i}(1-\pi_{i}%
)/E[\pi_{i}(1-\pi_{i})].$

It is straightforward to construct a DCDR estimator of $\beta_{0}$. Let
$\gamma_{0}(x_{i})=E[y_{i}|x_{i}]$ and $\alpha_{0}(x_{i})=-E[a_{i}|x_{i}]$ as
before, except that $a_{i}$ may now be a vector. Also let $I_{\ell}$ denote
the index set for the $\ell^{th}$ group, and $\hat{I}_{\ell}$ and $\tilde
{I}_{\ell}$ the index sets for the observations used to obtain $\hat{\gamma
}_{\ell}$ and $\tilde{\alpha}_{\ell}$ respectively. For any function $g(z)$
let
\[
\bar{F}\{g(z)\}=\frac{1}{\bar{n}_{\ell}}\sum_{i\in I_{\ell}}g(z_{i}),\hat
{F}\{g(z)\}=\frac{1}{\hat{n}_{\ell}}\sum_{i\in\hat{I}_{\ell}}g(z_{i}%
),\tilde{F}\{g(z)\}=\frac{1}{\tilde{n}_{\ell}}\sum_{i\in\tilde{I}_{\ell}%
}g(z_{i}).
\]
These represent sample averages over each of the groups of observations. Let
$\hat{\gamma}_{\ell}(x),$ $\hat{\alpha}_{\ell}(x),$ and $\tilde{\alpha}_{\ell
}(x)$ be series estimators of $\gamma_{0}(x)$ and $\alpha_{0}(x)$ given by
\begin{align*}
\hat{\gamma}_{\ell}(x)  &  =p(x)^{T}\hat{\delta}_{\ell},\hat{\alpha}_{\ell
}(x)=p(x)^{T}\hat{\delta}_{\ell\alpha},\tilde{\alpha}_{i}(x)=p(x)^{T}%
\tilde{\delta}_{\ell\alpha},\\
\hat{\delta}  &  =\hat{\Sigma}^{-1}\hat{h},\text{ }\hat{\delta}_{\alpha}%
=\hat{\Sigma}^{-}\hat{h}_{\alpha},\text{ }\tilde{\delta}_{\alpha}%
=\tilde{\Sigma}^{-}\tilde{h}_{\alpha},\text{ }\hat{\Sigma}=\hat{F}%
\{p(x)p(x)^{T}\},\text{ }\tilde{\Sigma}=\tilde{F}\{p(x)p(x)^{T}\},\\
\hat{h}  &  =\hat{F}\{p(x)y\},\text{ }\hat{h}_{\alpha}=\hat{F}\{p(x)a\},\text{
}\tilde{h}_{\alpha}=\tilde{F}\{p(x)a\}.
\end{align*}
The estimator we consider is%
\begin{equation}
\tilde{\beta}=\left(  \sum_{\ell=1}^{L}\sum_{i\in I_{\ell}}[a_{i}%
-\tilde{\alpha}_{\ell}(x_{i})][a_{i}-\hat{\alpha}_{\ell}(x_{i})]^{T}\right)
^{-1}\sum_{\ell=1}^{L}\sum_{i\in I_{\ell}}[a_{i}-\tilde{\alpha}_{\ell}%
(x_{i})][y_{i}-\hat{\gamma}_{\ell}(x_{i})]. \label{drpartlin}%
\end{equation}
This estimator can be thought of as an instrumental variables estimator with
left hand sides variable $y_{i}-\hat{g}_{i}(x_{i})$, right hand side variables
$a_{i}-\hat{\alpha}_{i}(x_{i}),$ and instruments $a_{i}-\tilde{\alpha}%
_{i}(x_{i}).$ Here the instrumental variables form is used to implement the
cross-fitting and not to correct for endogeneity. This form means that every
element of the matrix that is inverted and of the vector it is multiplying is
a DCDR estimator of an expected conditional covariance like that described earlier.

\bigskip

\textsc{Theorem 5: }\textit{If Assumptions 1 - 3 and 5 are satisfied,
}$\lambda_{0}(x)=E[y_{i}-a_{i}^{T}\beta_{0}|x_{i}=x]$\textit{ is Holder of
order }$s_{\gamma}$\textit{ and each component of }$E[a_{i}|x_{i}=x]$\textit{
is Holder of order }$s_{\alpha}$\textit{, }$H=E[Var(a_{i}|x_{i})]$\textit{
exists and is nonsingular, and }$\Omega=E[\{a_{i}-\alpha_{0}(x_{i}%
)\}\{a_{i}-\alpha_{0}(x_{i})\}^{T}\varepsilon_{i}^{2}]$\textit{ exists then
for }$\varepsilon_{i}=y_{i}-a_{i}^{T}\beta_{0}-\lambda_{0}(x_{i})$\textit{ and
}$\psi(z_{i})=H^{-1}(a_{i}-E[a_{i}|x_{i}])\varepsilon_{i},$%
\[
\sqrt{n}(\tilde{\beta}-\beta_{0})=\frac{1}{\sqrt{n}}\sum_{i=1}^{n}\psi
(z_{i})+O_{p}(\bar{\Delta}_{n}^{\ast}).
\]

\bigskip

The regularity conditions here are somewhat stronger than those of Donald and
Newey (1994), who do not require any restrictions on the marginal distribution
of $x_{i}$ nor use any sample splitting. This strengthening is useful to
achieve the fast remainder for partially linear projections rather than for
the coefficients $\beta_{0}$ in the conditional mean model $E[y_{i}%
|a_{i},x_{i}]=a_{i}^{T}\beta_{0}+\lambda_{0}(x_{i})$ of Donald in Newey
(1994). The upper bound on the rate at which $K$ can grow is slightly stricter
than in Donald and Newey (1994) due to the presence of the $\ln(K)$ term in
Assumption 5. Thus, under somewhat stronger conditions than those of Donald
and Newey (1994) the DCDR estimator of a partially linear projection has a
fast remainder just as in Donald and Newey (1994). Consequently, the estimator
will be root-n consistent under minimal conditions.

When the Robins et al. (2009) minimal condition $(s_{\gamma}+s_{\alpha
})/r>1/2$ holds, consider a spline with $\kappa>\max\{s_{\gamma},s_{\alpha
}\}-1$, so that $\zeta_{\gamma}+\zeta_{\alpha}=(s_{\gamma}+s_{\alpha})/r>1/2$.
Then there will exist a $K$ such that $\bar{\Delta}_{n}^{\ast}\longrightarrow
0$ and hence $\tilde{\beta}$ will be semiparametric efficient. Thus we see
that the DCDR estimator $\tilde{\beta}$ of equation (\ref{drpartlin}%
)\thinspace will be semiparametric efficient under nearly minimal conditions
and has a fast remainder term.

\bigskip\ \ 

\section{The Doubly Robust Estimator}

In this Section we show that the DCDR estimator has improved properties
relative to the plug-in estimator, in the sense that the remainder bounds are
smaller for the DCDR robust estimator. We have not yet been able to obtain the
fast remainder for the doubly robust estimator for general splines, for the
same reasons as for plug-in estimators.

Before giving results for series estimators we give a result that applies to
any doubly robust estimator of a linear functional. Let $\mathcal{A}_{n}$
denote an event that occurs with probability approaching one. For example,
$\mathcal{A}_{n}$ could include the set of data points where $\hat{\Sigma
}_{\ell}$ is nonsingular.

\bigskip

\textsc{Lemma 6:} \textit{If Assumptions 1 and 2 are satisfied, }$\hat{\gamma
}_{\ell}(x)$\textit{ and }$\hat{\alpha}_{\ell}(x)$\textit{ do not use
observations in }$I_{\ell}$\textit{, }$Var(y_{i}|x_{i})$\textit{ is bounded,
and there are }$\Delta_{n}^{m},$ $\Delta_{m}^{\gamma}$, and $\Delta
_{m}^{\alpha}$, \textit{such that for each }$\left(  \ell=1,...,L\right)
$,\textit{ }%
\begin{align*}
1(\mathcal{A}_{n})\left\{  \int[m(z,\hat{\gamma}_{\ell})-m(z,\gamma_{0}%
)]^{2}F_{0}(dz)\right\}  ^{1/2}  &  =O_{p}(\Delta_{n}^{m}),\\
1(\mathcal{A}_{n})\left\{  \int\alpha_{0}(x)^{2}[\hat{\gamma}_{\ell}%
(x)-\gamma_{0}(x)]^{2}F_{0}(dz)\right\}  ^{1/2}  &  =O_{p}(\Delta_{n}^{\gamma
}),\\
1(\mathcal{A}_{n})\left\{  \int[\tilde{\alpha}_{\ell}(x)-\alpha_{0}%
(x)]^{2}F_{0}(dz)\right\}  ^{1/2}  &  =O_{p}(\Delta_{n}^{\alpha}),
\end{align*}
\textit{ then}%
\[
\sqrt{n}(\tilde{\beta}-\beta_{0})=\frac{1}{\sqrt{n}}\sum_{i=1}^{n}\psi
(z_{i})-\frac{1}{\sqrt{n}}\sum_{\ell=1}^{L}\sum_{i\in I_{\ell}}[\tilde{\alpha
}_{\ell}(x_{i})-\alpha_{0}(x_{i})][\hat{\gamma}_{\ell}(x_{i})-\gamma_{0}%
(x_{i})]+O_{p}(\Delta_{n}^{m}+\Delta_{n}^{\gamma}+\Delta_{n}^{\alpha}).
\]

\bigskip

This result does not require that $\hat{\gamma}_{\ell}(x)$ and\textit{ }%
$\hat{\alpha}_{\ell}(x)$ be computed from different samples. It only uses the
sample splitting in averaging over different observations that are used to
construct $\hat{\gamma}_{\ell}$ and $\tilde{\alpha}_{\ell}.$ Also, it is known
from Newey, Hsieh, and Robins (1998, 2004) and Chernozhukov et. al. (2016)
that adding the adjustment term to the plug-in estimator makes the remainder
second order. The conclusion of Lemma 6 gives an explicit form of that result.
Under weak conditions that only involve mean-square convergence the doubly
robust estimator has a remainder that is the sum of three stochastic
equicontinuity remainders and the quadratic, split sample remainder involving
the product of the estimation remainders for the two nonparametric estimators
$\hat{\gamma}$ and $\tilde{\alpha}$.

For series estimators the DCDR estimator will have $\bar{\Delta}_{n}^{\ast}$
as its primary remainder for the Haar basis

\bigskip

\textsc{Theorem 7:} \textit{If Assumptions 1-6 are satisfied, }$\kappa
=0$\textit{, and }$K[\ln(n)]^{2}/n\longrightarrow0$ \textit{then}%
\[
\sqrt{n}(\tilde{\beta}-\beta_{0})=\frac{1}{\sqrt{n}}\sum_{i=1}^{n}\psi
(z_{i})+O_{p}(\bar{\Delta}_{n}^{\ast}+\Delta_{n}^{m}).
\]
\textit{If in addition }$d_{K}$\textit{ is bounded as a function of }%
$K$\textit{ and }$\zeta_{m}=\zeta_{\gamma}$\textit{ then }$\Delta_{n}^{m}\leq
C\bar{\Delta}_{n}^{\ast}$\textit{.}

\bigskip

One improvement of the DCDR estimator over the plug-in estimator is that the
remainder no longer contains the $K^{-\zeta_{\gamma}}\ln(n)$ term. The
elimination of this term is the direct result of the DCDR estimator having a
smaller remainder than the plug-in estimator.

For splines of order $\kappa>0$ we can obtain a result for the DCDR estimator
that improves on the plug-in remainder bound.

\bigskip

\textsc{Theorem 8:} \textit{If Assumptions 1-6 are satisfied} \textit{then}%
\[
\sqrt{n}(\tilde{\beta}-\beta_{0})=\frac{1}{\sqrt{n}}\sum_{i=1}^{n}\psi
(z_{i})+O_{p}(\bar{\Delta}_{n}^{\ast}+\Delta_{n}^{m}+\tilde{\Delta}%
_{n}),\tilde{\Delta}_{n}=\sqrt{\frac{K^{3}\left[  \ln(K)\right]  ^{2}%
(1+d_{K})}{n^{3}}}K^{(1/2)-\zeta_{\gamma}}.
\]
\textit{If in addition }$d_{K}$\textit{ is bounded as a function of }%
$K$\textit{ and }$\zeta_{m}=\zeta_{\gamma}$\textit{ then }$\Delta_{n}^{m}\leq
C\bar{\Delta}_{n}^{\ast}$\textit{.}

\bigskip

Here we see that the remainder bound for the DCDR estimator will generally be
smaller than the remainder bound for the plug-in estimator because the term
$K\ln(K)/n$ is raised to the 3/2 power rather than the 1/2 power. Here it
turns out that there will exist a $K$ such that all of the remainder terms go
to zero if%
\[
4\zeta_{\gamma}+3\zeta_{\alpha}\geq2.
\]
For example, if $s_{\gamma}=s_{\alpha}$ and $\kappa\geq\max\{s_{\gamma
},s_{\alpha}\}-1,$ this requires $s_{\gamma}>2r/7,$ which is only slightly
stronger than the $s_{\gamma}>r/4$ condition of Robins et al.(2009) that is
required for existence of a semiparametric efficient estimator. Also,
existence of $K$ such that the remainder will be of size no larger than
$\bar{\Delta}_{n}^{\ast}$ requires%
\[
2\zeta_{\gamma}+\zeta_{\alpha}\geq1.
\]
For example, if $\zeta_{\gamma}=\zeta_{\alpha}$ this requires $\zeta_{\gamma
}>1/3,$ which is weaker than the condition $\zeta_{\gamma}>1/2$ for the
remainder for the plug-in estimator. In these ways the DCDR estimator improves
on the plug-in estimator.

\bigskip

\section{Appendix}

This Appendix gives the proofs of the results in the body of the paper. We
begin with the proofs of Lemma 1 and Lemma 6 because they are not restricted
to series estimators.

\bigskip

\textbf{Proof of Lemma 1:} Define $\hat{\Delta}_{i\ell}=m(z_{i},\hat{\gamma
})-m(z_{i},\gamma_{0})-\bar{m}(\hat{\gamma}_{\ell})+\beta_{0}$ for $i\in
I_{\ell}$ and let $Z(I_{\ell})^{c}$ denote the set of observations $z_{i}$ for
$i\notin I_{\ell}$. Note that $E[\hat{\Delta}_{i\ell}|Z(I_{\ell})^{c}]=0$ by
construction for $i\in I_{\ell}.$ Also by independence of the observations,
$E[\hat{\Delta}_{i\ell}\hat{\Delta}_{j\ell}|Z(I_{\ell})^{c}]=0$ for $i,j\in
I_{\ell}.$ Furthermore, $E[\hat{\Delta}_{i\ell}^{2}|Z(I_{\ell})^{c}]\leq
\int[m(z,\hat{\gamma}_{\ell})-m(z,\gamma_{0})]^{2}F_{0}(dz)=O_{p}((\Delta
_{n}^{m})^{2})$ for $i\in I_{\ell}$. Then we have
\[
E[\left(  \frac{1}{\sqrt{n}}\sum_{i\in I_{\ell}}\hat{\Delta}_{i\ell}\right)
^{2}|Z(I_{\ell})^{c}]=\frac{1}{n}E[\left(  \sum_{i\in I_{\ell}}\hat{\Delta
}_{i\ell}\right)  ^{2}|Z(I_{\ell})^{c}]=\frac{\bar{n}_{\ell}}{n}E[\hat{\Delta
}_{i\ell}^{2}|Z(I_{\ell})^{c}]=O_{p}((\Delta_{n}^{m})^{2}).
\]
Therefore, by the Markov inequality we have $\sum_{i\in I_{\ell}}\hat{\Delta
}_{i\ell}/\sqrt{n}=O_{p}(\Delta_{n}^{m}).$ The first conclusion then follows from%

\[
\sqrt{n}(\hat{\beta}-\beta_{0})=\sum_{\ell=1}^{L}\frac{1}{\sqrt{n}}\sum_{i\in
I_{\ell}}\hat{\Delta}_{i\ell}+\frac{1}{\sqrt{n}}\sum_{i=1}^{n}[m(z_{i}%
,\gamma_{0})-\beta_{0}]+\sqrt{n}\sum_{\ell=1}^{L}\frac{\bar{n}_{\ell}}{n}%
[\bar{m}(\hat{\gamma}_{\ell})-\beta_{0}].
\]
For the second conclusion note by the subsamples being as close to equal size
as possible,%
\[
\frac{\bar{n}_{\ell}}{\hat{n}_{\ell}}=\frac{\bar{n}_{\ell}/n}{\hat{n}_{\ell
}/n}=\frac{1/L}{(L-1)/L}+O(n^{-1})=\frac{1}{(L-1)}+O(n^{-1}).
\]
Then by%
\begin{align*}
\sqrt{n}\sum_{\ell=1}^{L}\frac{\bar{n}_{\ell}}{n}[\bar{m}(\hat{\gamma}_{\ell
})-\beta_{0}]  &  =\frac{1}{\sqrt{n}}\sum_{\ell=1}^{L}\bar{n}_{\ell}%
\sqrt{\frac{1}{\hat{n}_{\ell}}}\sqrt{\hat{n}_{\ell}}[\bar{m}(\hat{\gamma
}_{\ell})-\beta_{0}]=\sum_{\ell=1}^{L}\frac{\bar{n}_{\ell}}{\hat{n}_{\ell}%
}\frac{1}{\sqrt{n}}\sum_{i\notin I_{\ell}}\phi(z_{i})+O_{p}(\Delta_{n}^{\phi
})\\
&  =\frac{1}{L-1}\frac{1}{\sqrt{n}}\sum_{\ell=1}^{L}\sum_{i\notin I_{\ell}%
}\phi(z_{i})+O_{p}(\Delta_{n}^{\phi}+n^{-1})\\
&  =\frac{1}{L-1}\frac{1}{\sqrt{n}}\sum_{\ell=1}^{L}(\sum_{i=1}^{n}\phi
(z_{i})-\sum_{i\in I_{\ell}}\phi(z_{i}))+O_{p}(\Delta_{n}^{\phi}+n^{-1})\\
&  =\frac{1}{\sqrt{n}}\sum_{i=1}^{n}\phi(z_{i})+O_{p}(\Delta_{n}^{\phi}%
+n^{-1}).
\end{align*}

The conclusion then follows by the triangle inequality. \textit{Q.E.D}.

\bigskip

\textbf{Proof of Lemma 6:} By adding and subtracting terms it follows that for
$\varepsilon_{i}=y_{i}-\gamma_{0}(x_{i})$ and $\phi(z_{i})=\alpha_{0}%
(x_{i})[y_{i}-\gamma_{0}(x_{i})]$%
\begin{align*}
\tilde{\alpha}_{\ell}(x_{i})[y_{i}-\hat{\gamma}_{\ell}(x_{i})]  &  =\phi
(z_{i})-\alpha_{0}(x_{i})[\hat{\gamma}(x_{i})-\gamma_{0}(x_{i})]+[\tilde
{\alpha}_{\ell}(x_{i})-\alpha_{0}(x_{i})]\varepsilon_{i}\\
&  -[\tilde{\alpha}_{\ell}(x_{i})-\alpha_{0}(x_{i})][\hat{\gamma}%
(x_{i})-\gamma_{0}(x_{i})].
\end{align*}
The first conclusion of Lemma 1 with $m(z,\gamma)=\alpha_{0}(x)\gamma(x)$
gives%
\[
\frac{1}{\sqrt{n}}\sum_{\ell=1}^{L}\sum_{i\in I_{\ell}}\alpha_{0}(x_{i}%
)[\hat{\gamma}_{\ell}(x_{i})-\gamma_{0}(x_{i})]=\sqrt{n}\sum_{\ell=1}^{L}%
\frac{\bar{n}_{\ell}}{n}\int\alpha(x)[\hat{\gamma}_{\ell}(x)-\gamma
_{0}(x)]F_{0}(dx)+O_{p}(\Delta_{n}^{\gamma}).
\]
Assumption 1 and the first conclusion of Lemma 1 also give%
\begin{align*}
\sqrt{n}\sum_{\ell=1}^{L}\frac{\bar{n}_{\ell}}{n}\int\alpha(x)[\hat{\gamma
}_{\ell}(x)-\gamma_{0}(x)]F_{0}(dx)  &  =\sqrt{n}\sum_{\ell=1}^{L}\frac
{\bar{n}_{\ell}}{n}[\bar{m}(\hat{\gamma}_{\ell})-\beta_{0}]\\
&  =\frac{1}{\sqrt{n}}\sum_{\ell=1}^{L}\sum_{i\in I_{\ell}}[m(z_{i}%
,\hat{\gamma}_{\ell})-m(z_{i},\gamma_{0})]+O_{p}(\Delta_{n}^{m}).
\end{align*}
In addition, if we take $\gamma=\alpha$ and $m(z,\alpha)=\alpha(x)\varepsilon$
then $\int m(z,\alpha)F_{0}(dz)=0$, so that by Lemma 1,%
\[
\frac{1}{\sqrt{n}}\sum_{\ell=1}^{L}\sum_{i\in I_{\ell}}[\tilde{\alpha}_{\ell
}(x_{i})-\alpha_{0}(x_{i})]\varepsilon_{i}=O_{p}(\Delta_{n}^{\alpha}).
\]
Then collecting terms we have%
\begin{align*}
\sqrt{n}(\tilde{\beta}-\beta_{0})  &  =\frac{1}{\sqrt{n}}\sum_{i=1}%
^{n}[m(z_{i},\gamma_{0})-\beta_{0}]\\
&  +\frac{1}{\sqrt{n}}\sum_{\ell=1}^{L}\sum_{i\in I_{\ell}}\{m(z_{i}%
,\hat{\gamma})-m(z_{i},\gamma_{0})+\tilde{\alpha}_{\ell}(x_{i})[y_{i}%
-\hat{\gamma}_{\ell}(x_{i})]\}\\
&  =\frac{1}{\sqrt{n}}\sum_{i=1}^{n}\psi(z_{i})+\frac{1}{\sqrt{n}}\sum
_{\ell=1}^{L}\sum_{i\in I_{\ell}}\alpha_{0}(x_{i})[\hat{\gamma}_{\ell}%
(x_{i})-\gamma_{0}(x_{i})]+O_{p}(\Delta_{n}^{m}+\Delta_{n}^{\gamma})\\
&  \frac{1}{\sqrt{n}}\sum_{\ell=1}^{L}\sum_{i\in I_{\ell}}\{-\alpha_{0}%
(x_{i})[\hat{\gamma}(x_{i})-\gamma_{0}(x_{i})]+[\tilde{\alpha}_{\ell}%
(x_{i})-\alpha_{0}(x_{i})]\varepsilon_{i}\}\\
&  -\frac{1}{\sqrt{n}}\sum_{\ell=1}^{L}\sum_{i\in I_{\ell}}[\tilde{\alpha
}_{\ell}(x_{i})-\alpha_{0}(x_{i})][\hat{\gamma}(x_{i})-\gamma_{0}(x_{i})]\\
&  =\frac{1}{\sqrt{n}}\sum_{i=1}^{n}\psi(z_{i})-\frac{1}{\sqrt{n}}\sum
_{\ell=1}^{L}\sum_{i\in I_{\ell}}[\tilde{\alpha}_{\ell}(x_{i})-\alpha
_{0}(x_{i})][\hat{\gamma}(x_{i})-\gamma_{0}(x_{i})]\\
&  +O_{p}(\Delta_{n}^{m}+\Delta_{n}^{\gamma}+\Delta_{n}^{\alpha}).Q.E.D.
\end{align*}

\bigskip

We now turn to proofs of the results involving series estimators. Let
$\Sigma=E[p(x_{i})p(x_{i})^{T}]$. It follows from Assumption 3 that $\Sigma$
is nonsingular, so we can replace $p(x)$ by $\Sigma^{-1/2}p(x)$ and so
normalize $\Sigma=I$ without changing the assumptions. We impose this
normalization throughout. Also, throughout the Appendix $C$ will denote a
generic constant not depending on $n$ or $K.$

We next prove the key result in eq. (\ref{key result}) for a zero order
spline. Let $r(x)=\gamma_{0}(x)-\gamma_{K}(x)$ and $\hat{h}_{2}=\sum_{i=1}%
^{n}p(x_{i})r(x_{i})/n$ as in the body of the paper. Also let $\left\Vert
A\right\Vert _{op}$ denote the operator norm of a symmetric matrix $A$, being
the largest absolute value of eigenvalues.

\bigskip

\textsc{Lemma A1: }\textit{If Assumptions 1-6 are satisfied, }$\kappa=0,$
$K[\ln(n)]^{2}/n\longrightarrow0$, \textit{then for }$\hat{U}=\sum_{j=0}%
^{J-1}(I-\hat{\Sigma})^{j}\hat{h}_{2},$ $\hat{W}=\hat{\Sigma}^{-1}%
(I-\hat{\Sigma})^{J}\hat{h}_{2},$ $J=int[\ln(n)]$ \textit{and any constant
}$\Delta>0,$%
\[
\left\Vert E[\hat{U}\hat{U}^{T}]\right\Vert _{op}\leq C\frac{K^{-2\zeta
_{\gamma}}[\ln(n)]^{2}}{n},\hat{W}^{T}\hat{W}=o_{p}(n^{-\Delta}).
\]

Proof: Let $Q_{i}=p(x_{i})p(x_{i})^{T},$ $\Delta_{i}=I-Q_{i}$, and
$h_{i}=p(x_{i})r(x_{i})$. Note that $E[\Delta_{i}]=0$ and $E[h_{i}]=0.$ For
each $j$ let $L=2j+2$. Let $\hat{U}_{j}=(I-\hat{\Sigma})^{j}\hat{h}_{2}.$ Then
we have%
\[
E[\hat{U}_{j}\hat{U}_{j}^{T}]=\frac{1}{n^{2j+2}}\sum_{i_{1},...,i_{L}=1}%
^{n}E[\left(  \Pi_{\ell=1}^{j}\Delta_{i_{\ell}}\right)  h_{i_{j+1}}h_{i_{j+2}%
}^{T}\left(  \Pi_{\ell=j+3}^{L}\Delta_{i_{\ell}}\right)  ].
\]
Consider any $(i_{1},...,i_{L})$ such that $i_{j+1}\neq i_{j+2}.$ Let
$i^{\ast}=i_{j+1}$ and let $Z_{i^{\ast}}^{c}$ denote the vector of
observations other than $z_{i^{\ast}}.$ Note that%
\[
E[\left(  \Pi_{\ell=1}^{j}\Delta_{i_{\ell}}\right)  h_{i_{j+1}}h_{i_{j+2}}%
^{T}\left(  \Pi_{\ell=j+3}^{L}\Delta_{i_{\ell}}\right)  ]=E[E[\left(
\Pi_{\ell=1}^{j}\Delta_{i_{\ell}}\right)  h_{i^{\ast}}h_{i_{j+2}}^{T}\left(
\Pi_{\ell=j+3}^{L}\Delta_{i_{\ell}}\right)  |Z_{i^{\ast}}^{c}]].
\]
We proceed to show that%
\[
E[\left(  \Pi_{\ell=1}^{j}\Delta_{i_{\ell}}\right)  h_{i^{\ast}}h_{i_{j+2}%
}^{T}\left(  \Pi_{\ell=j+3}^{L}\Delta_{i_{\ell}}\right)  |Z_{i^{\ast}}%
^{c}]=0.
\]
Note that conditional on $Z_{i^{\ast}}^{c}$ we can treat all terms where
$i_{\ell}\neq i^{\ast}$ as constant. Also, because $i_{j+1}\neq i_{j+2}$ all
terms where $i_{\ell}=i^{\ast}$ depend only on $p(x_{i^{\ast}}).$ Therefore
for the scalar $r(x)=\gamma_{0}(x)-\gamma_{K}(x)$ we have%
\[
E[\left(  \Pi_{\ell=1}^{j}\Delta_{i_{\ell}}\right)  h_{i^{\ast}}h_{i_{j+2}%
}^{T}\left(  \Pi_{\ell=j+3}^{L}\Delta_{i_{\ell}}\right)  |Z_{i^{\ast}}%
^{c}]=E[A_{1}(p(x_{i^{\ast}}))p(x_{i^{\ast}})r(x_{i^{\ast}})A_{2}%
(p(x_{i^{\ast}}))]=E[A(p(x_{i^{\ast}}))r(x_{i^{\ast}})],
\]
where $A_{1}(p)$ and $A_{2}(p)$ are $K\times K$ and $1\times K$ matrices of
functions of $p$ and $A(p)=A_{1}(p)pA_{2}(p).$ Let $X_{k}$ denote the interval
where $p_{k}(x)$ is nonzero. Note that $p_{k}(x)=1(x\in X_{k})c_{k}$ for a
constant $c_{k}$, and hence%
\[
A(p(x_{i^{\ast}}))=\sum_{k=1}^{K}A_{k}1(x_{i^{\ast}}\in X_{k}),A_{k}%
=A((0,...,0,c_{k},0,...,0)^{T}).
\]
Therefore by orthogonality of each $p_{k}(x_{i})$ with $r(x_{i})$ in the
population,%
\[
E[\left(  \Pi_{\ell=1}^{j}\Delta_{i_{\ell}}\right)  h_{i^{\ast}}h_{i_{j+2}%
}^{T}\left(  \Pi_{\ell=j+3}^{L}\Delta_{i_{\ell}}\right)  |Z_{i^{\ast}}%
^{c}]=\sum_{k=1}^{K}A_{k}E[1(x_{i^{\ast}}\in X_{k})r(x_{i^{\ast}})]=\sum
_{k=1}^{K}A_{k}c_{k}^{-1}E[p_{k}(x_{i^{\ast}})r(x_{i^{\ast}})]=0.
\]
Therefore by iterated expectations, if $i_{j+1}\neq i_{j+2}$ we have%
\[
E[\left(  \Pi_{\ell=1}^{j}\Delta_{i_{\ell}}\right)  h_{i_{j+1}}h_{i_{j+2}}%
^{T}\left(  \Pi_{\ell=j+3}^{L}\Delta_{i_{\ell}}\right)  ]=0.
\]
It then follows that for $\Psi=E[h_{i_{j+1}}h_{i_{j+1}}^{T}]=E[r(x_{i}%
)^{2}p(x_{i})p(x_{i})^{T}]$ and $\tilde{\Delta}_{i_{j+1}}=h_{i_{j+1}%
}h_{i_{j+1}}^{T}-\Psi,$%
\begin{align*}
E[\hat{U}_{j}\hat{U}_{j}^{T}]  &  =\frac{1}{n^{2j+2}}\sum_{i_{1}%
,..,i_{j+1},i_{j+3}....,i_{L}=1}^{n}E[\left(  \Pi_{\ell=1}^{j}\Delta_{i_{\ell
}}\right)  h_{i_{j+1}}h_{i_{j+1}}^{T}\left(  \Pi_{\ell=j+3}^{L}\Delta
_{i_{\ell}}\right)  ]=T_{1}^{j}+T_{2}^{j},\\
T_{1}^{j}  &  =\frac{1}{n^{2j+1}}\sum_{i_{1},..,i_{j},i_{j+3}....,i_{L}=1}%
^{n}E[\left(  \Pi_{\ell=1}^{j}\Delta_{i_{\ell}}\right)  \Psi\left(  \Pi
_{\ell=j+3}^{L}\Delta_{i_{\ell}}\right)  ],\\
T_{2}^{j}  &  =\frac{1}{n^{2j+2}}\sum_{i_{1},..,i_{j+1},i_{j+3}....,i_{L}%
=1}^{n}E[\left(  \Pi_{\ell=1}^{j}\Delta_{i_{\ell}}\right)  \tilde{\Delta
}_{i_{j+1}}\left(  \Pi_{\ell=j+3}^{L}\Delta_{i_{\ell}}\right)  ].
\end{align*}

Consider first $T_{2}^{j}$. Note that $\Delta_{i}$ and $\tilde{\Delta}_{i}$
are diagonal matrices, so that $E[\left(  \Pi_{\ell=1}^{j}\Delta_{i_{\ell}%
}\right)  \tilde{\Delta}_{i_{j+1}}\left(  \Pi_{\ell=j+3}^{L}\Delta_{i_{\ell}%
}\right)  ]$ is a diagonal matrix, with $k^{th}$ diagonal element given by
$E[\left(  \Pi_{\ell=1}^{j}\Delta_{k,i_{\ell}}\right)  \tilde{\Delta
}_{k,i_{j+1}}\left(  \Pi_{\ell=j+3}^{L}\Delta_{k,i_{\ell}}\right)  ]$ where
\[
\Delta_{k,i}=p_{k}(x_{i})^{2}-E[p_{k}(x_{i})^{2}],\tilde{\Delta}_{k,i_{j+1}%
}=r(x_{i})^{2}p_{k}(x_{i})^{2}-E[r(x_{i})^{2}p_{k}(x_{i})^{2}].
\]
The largest absolute value of the eigenvalues of a diagonal matrix is the
maximum of the absolute values of the diagonal elements, so it suffices to
show that the conclusion holds for these diagonal elements. We will consider
the $k^{th}$ diagonal element but for notational convenience drop the $k$
subscript in what follows.

Note that $p_{k}(x_{i})^{2}\leq BK$ for some $B$ that does not vary with $k$
or $j.$ Also, for any random variable $Y_{i}$ and $\mu=E[Y_{i}],$ note that by
Jensen's inequality, $\left\vert \mu\right\vert ^{s}\leq E[|Y_{i}|^{s}]$ for
$s\geq1.$ Then for any positive $s$,
\[
E[|Y_{i}-\mu|^{s}]\leq E[\left(  |Y_{i}|+|\mu|\right)  ^{s}]\leq
E[2^{s-1}\left(  |Y_{i}|^{s}+|\mu|^{s}\right)  ]\leq2^{s-1}\left(
E[|Y_{i}|^{s}]+|\mu|^{s}\right)  ]\leq2^{s}E[|Y_{i}|^{s}]
\]
Then for any positive integer $s$, by the triangle inequality and the
definitions of $\Delta_{i},$%
\begin{equation}
\left\vert E[\Delta_{i}^{s}]\right\vert \leq2^{s}E[p_{k}(x_{i})^{2s}]\leq
2^{s}(BK)^{s-1}E[p_{k}(x_{i})^{2}]\leq(4BK)^{s-1}\leq(CK)^{s-1}.
\label{powbound1}%
\end{equation}
Also, by $r(x_{i})^{2}\leq DK^{-2\zeta_{\gamma}}$ we have%
\begin{align}
\left\vert E[(\Delta_{i})^{s}\tilde{\Delta}_{i}]\right\vert  &  \leq
E[\left\vert \Delta_{i}\right\vert ^{s}(r(x_{i})^{2}p_{k}(x_{i})^{2}%
+E[r(x_{i})^{2}p_{k}(x_{i})^{2}])]\label{powbound}\\
&  \leq E[(p_{k}(x_{i})^{2}+E[p_{k}(x_{i})^{2}])^{s+1}]DK^{-2\zeta_{\gamma}%
}\nonumber\\
&  \leq2^{s+1}E[p_{k}(x_{i})^{2s+2}]DK^{-2\zeta_{\gamma}}\leq2^{s+1}%
(BK)^{s}DK^{-2\zeta_{\gamma}}\nonumber\\
&  \leq(4(D+1)BK)^{s}K^{-2\zeta_{\gamma}}\leq(CK)^{s}K^{-2\zeta_{\gamma}%
}.\nonumber
\end{align}
\qquad\qquad\qquad

Next consider%
\[
T_{2}^{j}=\frac{1}{n^{2j+2}}\sum_{i_{1},..,i_{j+1},i_{j+3}....,i_{2j+2}=1}%
^{n}E[\left(  \Pi_{\ell=1}^{j}\Delta_{i_{\ell}}\right)  \tilde{\Delta
}_{i_{j+1}}\left(  \Pi_{\ell=j+3}^{L}\Delta_{i_{\ell}}\right)  ].
\]
The only terms in this sum that are nonzero are those where every index
$i_{\ell}$ is equal to at least one other index $i_{\ell^{\prime}}$, i.e.
where each index is "matched" with at least one other. Let $\tilde{\imath
}=(i_{1},..,i_{j+1},i_{j+3}....,i_{2j+2})^{T}$ denote the $2j+1$ dimensional
vector of indices where each $i_{\ell}$ is an integer in $[1,n].$ Let
$\Upsilon_{d}$ denote a set of all such $\tilde{\imath}$ with specified
indices that are equal to each other, but those matched indices are not equal
to any other indices. For example, one $\Upsilon_{d}$ is the set of
$\tilde{\imath}$ with $i_{1}=i_{j+1}=i_{j+3}=\cdots=i_{2J+2}$ and another is
the set of $\tilde{\imath}$ with $i_{1}=i_{2},i_{3}=\cdots=i_{2J+2},i_{2}\neq
i_{3}.$ For each $d$ each group of index coordinates that are equal to each
other can be thought of as a group of matching indices that we index by
$g_{d}.$ Let $m_{g_{d}}$ denote the number of indices in group $g_{d}$ and
$G_{d}$ denote the total number of groups. Note that the total number of
indices is $2j+1=\sum_{g_{d}=1}^{G_{d}}m_{g_{d}}$. Also, by eqs.
(\ref{powbound1}) and (\ref{powbound}) for each $\tilde{\imath}\in \Upsilon
_{d}$ we have
\[
|E[\left(  \Pi_{\ell=1}^{j}\Delta_{i_{\ell}}\right)  \tilde{\Delta}_{i_{j+1}%
}\left(  \Pi_{\ell=j+3}^{L}\Delta_{i_{\ell}}\right)  ]|\leq K^{-2\zeta
_{\gamma}}%
{\displaystyle\prod\limits_{g_{d}=1}^{G_{d}}}
\left(  CK\right)  ^{m_{gd}-1}=K^{-2\zeta_{\gamma}}\left(  CK\right)
^{2j+1-G_{d}}.
\]
Also, the number of indices in $\Upsilon_{d}$ is less than or equal to
$n^{G_{d}}$ since each match can be regarded as a single index. Therefore,%
\begin{align*}
\left\vert \frac{1}{n^{2j+2}}\sum_{\tilde{\imath}\in \Upsilon_{d}}^{n}E[\left(
\Pi_{\ell=1}^{j}\Delta_{i_{\ell}}\right)  \tilde{\Delta}_{i_{j+1}}\left(
\Pi_{\ell=j+3}^{L}\Delta_{i_{\ell}}\right)  ]\right\vert  &  \leq\frac
{1}{n^{2j+2}}\sum_{\tilde{\imath}\in \Upsilon_{d}}^{n}\left\vert E[\left(
\Pi_{\ell=1}^{j}\Delta_{i_{\ell}}\right)  \tilde{\Delta}_{i_{j+1}}\left(
\Pi_{\ell=j+3}^{L}\Delta_{i_{\ell}}\right)  ]\right\vert \\
&  \leq\left(  \frac{1}{n^{2j+2}}\right)  n^{G_{d}}K^{-2\zeta_{\gamma}}\left(
CK\right)  ^{2j+1-G_{d}}\\
&  =\frac{1}{n}K^{-2\zeta_{\gamma}}\left(  \frac{CK}{n}\right)  ^{2j+1-G_{d}}.
\end{align*}
By hypothesis $K/n\longrightarrow0$ so that for large enough $n$ we have
$CK/n<1$. For such $n$ we have $\left(  CK/n\right)  ^{2j+1-G_{d}}$ decreasing
in $G_{d}.$ Also, the largest $G_{d}$ is $j$, because each group must contain
at least two elements. Therefore, for large enough $n$ we have%
\[
\left\vert \frac{1}{n^{2j+2}}\sum_{\tilde{\imath}\in \Upsilon_{d}}^{n}E[\left(
\Pi_{\ell=1}^{j}\Delta_{i_{\ell}}\right)  \tilde{\Delta}_{i_{j+1}}\left(
\Pi_{\ell=j+3}^{L}\Delta_{i_{\ell}}\right)  ]\right\vert \leq\frac{1}%
{n}K^{-2\zeta_{\gamma}}\left(  \frac{CK}{n}\right)  ^{j+1}.
\]
Note that the bound on the right does not depend on $d$. Let $D$ denote the
total number of possible $\Upsilon_{d}$. Then since $E[\left(  \Pi_{\ell
=1}^{j}\Delta_{i_{\ell}}\right)  \tilde{\Delta}_{i_{j+1}}\left(  \Pi
_{\ell=j+3}^{L}\Delta_{i_{\ell}}\right)  ]=0$ if $\tilde{\imath}\notin%
\cup_{d=1}^{D}\Upsilon_{d}$ we have
\[
\left\vert T_{2}^{j}\right\vert \leq\sum_{d=1}^{D}\left\vert \frac{1}%
{n^{2j+2}}\sum_{\tilde{\imath}\in \Upsilon_{d}}^{n}E[\left(  \Pi_{\ell=1}%
^{j}\Delta_{i_{\ell}}\right)  \tilde{\Delta}_{i_{j+1}}\left(  \Pi_{\ell
=j+3}^{L}\Delta_{i_{\ell}}\right)  ]\right\vert \leq\frac{D}{n}K^{-2\zeta
_{\gamma}}\left(  \frac{CK}{n}\right)  ^{j+1}.
\]
Note that there are exactly $j^{2j+1}$ ways of forming $2j+1$ indices into $j$
groups. Ignoring the fact that we can exclude ways where any group has only
one index we have the bound $D\leq j^{2j+1}$. Plugging in this bound into the
above inequality and maximizing over diagonal elements gives%
\[
\left\Vert T_{2}^{j}\right\Vert _{op}\leq\frac{j^{2j+1}K^{-2\zeta_{\gamma}}%
}{n}\left(  \frac{CK}{n}\right)  ^{j+1}.
\]
Arguing similarly for $T_{1}^{j}$ gives%
\[
\left\Vert T_{1}^{j}\right\Vert _{op}\leq\frac{j^{2j}K^{-2\zeta_{\gamma}}}%
{n}\left(  \frac{CK}{n}\right)  ^{j},
\]
where we take $0^{0}=1.$

Next note that by $K\ln(n)^{2}/n\longrightarrow0$ we have $CK/n\leq
1/[2\ln(n)^{2}]$ for large enough $n.$ Also, $j/\ln(n)\leq1$ for all $j<J.$
Then for $n$ large enough%
\[
\sum_{j=0}^{J-1}j^{2j}\left(  \frac{CK}{n}\right)  ^{j}\leq\sum_{j=0}%
^{J-1}j^{2j}\left(  \frac{1}{2\ln(n)^{2}}\right)  ^{j}\leq\sum_{j=0}%
^{J-1}\left(  \frac{j}{\ln(n)}\right)  ^{2j}\left(  \frac{1}{2}\right)
^{j}\leq\sum_{j=0}^{J-1}\left(  \frac{1}{2}\right)  ^{j}\leq\sum_{j=0}%
^{\infty}\left(  \frac{1}{2}\right)  ^{j}=\frac{1}{1-\varepsilon_{n}}\leq2.
\]
Similarly it follows that for large enough $n,$%
\[
\sum_{j=0}^{J-1}j^{2j+1}\left(  \frac{CK}{n}\right)  ^{j+1}\leq\frac{1}%
{2\ln(n)}\sum_{j=0}^{J-1}\left(  \frac{j}{\ln(n)}\right)  ^{2j+1}\left(
\frac{1}{2}\right)  ^{j}\leq\frac{1}{\ln(n)}.
\]
Then we have for large enough $n$,
\begin{align*}
\left\Vert \sum_{j=0}^{J-1}E[\hat{U}_{j}\hat{U}_{j}^{T}]\right\Vert _{op}  &
\leq\left\Vert \sum_{j=0}^{J-1}\left(  T_{1}^{j}+T_{2}^{j}\right)  \right\Vert
_{op}\leq\sum_{j=0}^{J-1}\left(  \left\Vert T_{1}^{j}\right\Vert
_{op}+\left\Vert T_{2}^{j}\right\Vert _{op}\right) \\
&  \leq\frac{K^{-2\zeta_{\gamma}}}{n}\left(  2+\frac{1}{\ln(n)}\right)
\leq\frac{CK^{-2\zeta_{\gamma}}}{n}.
\end{align*}

Also by the Cauchy Schwartz inequality, $\hat{U}\hat{U}^{T}=\left(  \sum
_{j=0}^{J-1}\hat{U}_{j}\right)  \left(  \sum_{j=0}^{J-1}\hat{U}_{j}\right)
^{T}\leq J^{2}\sum_{j=0}^{J-1}\hat{U}_{j}\hat{U}_{j}^{T}.$ Therefore, for
large enough $n,$%
\[
\left\Vert E[\hat{U}\hat{U}^{T}]\right\Vert _{op}\leq J^{2}\left\Vert
\sum_{j=0}^{J-1}E[\hat{U}_{j}\hat{U}_{j}^{T}]\right\Vert _{op}\leq\frac
{C\ln(n)^{2}K^{-2\zeta_{\gamma}}}{n},
\]
giving the first conclusion.

For the second conclusion note that for any $\Delta>0$,%
\[
\ln\{n^{\Delta}[\ln(n)]^{-2\ln(n)+2}\}=\ln(n)[\Delta-2\ln(\ln(n))]+2\ln
(\ln(n))\longrightarrow-\infty.
\]
It follows that $[\ln(n)]^{-2\ln(n)+2}=o(n^{-\Delta})$ for any $\Delta$. Also,
by $K/n=o(\left[  1/\ln(n)\right]  ^{2})$ we have $K\ln\left(  K\right)
/n=o(1/\ln(n)),$ so that%
\[
\left(  \frac{K\ln(K)}{n}\right)  ^{2J}=o([\ln(n)]^{-2int(\ln(n))}%
)=o([\ln(n)]^{-2(\ln(n))+2})=o(n^{-\Delta}),
\]
for any $\Delta>0$. Then we have%
\[
\hat{1}\hat{W}^{T}\hat{W}\leq4\hat{h}_{2}^{T}(I-\hat{\Sigma})^{2J}\hat{h}%
_{2}\leq4\hat{h}_{2}^{T}\hat{h}_{2}\left\Vert I-\hat{\Sigma}\right\Vert
_{op}^{2J}=O_{p}(\frac{K^{1-2\zeta_{\gamma}}}{n}\left[  \frac{K\ln(K)}%
{n}\right]  ^{2J})=o_{p}(n^{-\Delta}),
\]
for any $\Delta>0$ by Rudelson's (1999) law of large numbers for random
matrices, giving the second conclusion. $Q.E.D.$

\bigskip

In the Appendix we focus on one subset $\bar{I}=I_{\ell}$ of observations and
let $\hat{I}$ and $\tilde{I}$ denote the observations used to compute
$\hat{\delta}$ and $\tilde{\delta}_{\alpha}$ respectively. Let $\bar{n},$
$\hat{n},$ $\tilde{n}$ denote the number of elements of $\bar{I},$ $\hat{I}$,
and $\tilde{I}$ respectively and%
\[
\bar{F}\{g(z)\}=\frac{1}{\bar{n}}\sum_{i\in\bar{I}}g(z_{i}),\hat
{F}\{g(z)\}=\frac{1}{\hat{n}}\sum_{i\in\hat{I}}g(z_{i}),\tilde{F}%
\{g(z)\}=\frac{1}{\tilde{n}}\sum_{i\in\tilde{I}}g(z_{i}),
\]
denote averages over the respective subsets of observations.

Next we make a few definitions we will use throughout. Let $\zeta_{\gamma},$
$\zeta_{\alpha},$ $\delta,$ $\gamma_{K},$ $\delta_{\alpha},$ and $\alpha_{K}$
be as defined following Assumption 4. Also, let
\begin{align*}
\varepsilon_{i}  &  =y_{i}-\gamma_{0}(x_{i}),r_{i}=\gamma_{0}(x_{i}%
)-\gamma_{K}(x_{i}),\eta_{i}=v(z_{i})-p(x_{i})\alpha_{0}(x_{i}),r_{i}^{\alpha
}=\alpha_{0}(x_{i})-\alpha_{K}(x_{i}),\\
\hat{h}_{1}  &  =\hat{F}\{p(x)\varepsilon\},\hat{h}_{2}=\hat{F}%
\{p(x)r\},\tilde{h}_{1}^{\alpha}=\tilde{F}\{\eta\},\tilde{h}_{2}^{\alpha
}=\tilde{F}\{p(x)r^{\alpha}\},\hat{\Sigma}=\hat{F}\{p(x)p(x)^{T}%
\},\tilde{\Sigma}=\tilde{F}\{p(x)p(x)^{T}\},\\
\hat{\Delta}_{1}  &  =\hat{\Sigma}^{-}\hat{h}_{1},\hat{\Delta}_{2}=\hat
{\Sigma}^{-}\hat{h}_{2},\tilde{\Delta}_{1}^{\alpha}=\tilde{\Sigma}^{-}%
\tilde{h}_{1}^{\alpha},\tilde{\Delta}_{2}^{\alpha}=\tilde{\Sigma}^{-}\tilde
{h}_{2}^{\alpha},\bar{\Sigma}=\bar{F}\{p(x)p(x)^{T}\},
\end{align*}
One\ piece of algebra we will use throughout is that, when $\hat{\Sigma}$ and
$\tilde{\Sigma}$ are nonsingular, by adding and subtracting $\hat{\Sigma}%
^{-1}\hat{F}\{p(x)\gamma_{0}(x)\}$ and $\tilde{\Sigma}^{-1}\tilde
{F}\{p(x)\alpha_{0}(x_{i})\}$ respectively we have%
\begin{equation}
\hat{\delta}-\delta=\hat{\Delta}_{1}+\hat{\Delta}_{2},\tilde{\delta}_{\alpha
}-\delta_{\alpha}=\tilde{\Delta}_{1}^{\alpha}+\tilde{\Delta}_{2}^{\alpha}.
\label{coeff}%
\end{equation}

Some properties of these objects will be useful in the proofs to follow. We
collect these properties in the following result. Let $\hat{1}$ and $\tilde
{1}$ denote the indicator function that the smallest eigenvalue of
$\hat{\Sigma}$ or $\tilde{\Sigma}$ is larger than $1/2$ respectively. As in
Belloni et al.(2015) $\Pr(\hat{1}=1)\longrightarrow1$ and $\Pr(\bar
{1}=1)\longrightarrow1$. Also, let $\hat{Z}^{c},$ $\tilde{Z}^{c},$ $\bar
{Z}^{c}$ denote all the other observations other than those indexed by
$\hat{I},$ $\tilde{I},$ or $\bar{I}$ respectively and $X=(x_{1},...,x_{n})$.

\bigskip

\textsc{Lemma A2:} \textit{If Assumptions 1-6 are is satisfied then}%
\begin{align*}
\text{i) }\hat{1}\left\Vert \hat{\Delta}_{1}\right\Vert  &  =O_{p}\left(
\frac{K}{n}\right)  ;\text{ ii) }\hat{1}\left\Vert \hat{\Delta}_{2}\right\Vert
=o_{p}\left(  K^{-2\zeta_{\gamma}}\frac{K}{n}\right)  ;\text{ }\\
\text{iii) }\hat{1}\left\Vert \hat{\Delta}_{1}^{\alpha}\right\Vert  &
=O_{p}\left(  \frac{(1+d_{K})K}{n}\right)  ,\text{ iv) }\hat{1}\left\Vert
\hat{\Delta}_{2}^{\alpha}\right\Vert =O_{p}\left(  \frac{K}{n}\right)  ,\\
\text{v) }\hat{1}\left\Vert \hat{\delta}-\delta\right\Vert ^{2}  &
=O_{p}\left(  \frac{K}{n}\right)  \text{; vi) }\tilde{1}\left\Vert
\tilde{\delta}_{\alpha}-\delta_{\alpha}\right\Vert ^{2}=O_{p}\left(
\frac{d_{K}K}{n}\right)  ,\text{ }\\
\text{vii) }\hat{1}E[\hat{\Delta}_{1}\hat{\Delta}_{1}^{T}|X,\hat{Z}^{c}]  &
\leq\frac{C}{n}I,\text{ viii) }\hat{1}\int[\hat{\gamma}(x)-\gamma_{0}%
(x)]^{2}F_{0}(dx)=O_{p}(\frac{K}{n}+K^{-2\zeta_{\gamma}}),\text{ }\\
\text{ix) }\tilde{1}\int[\tilde{\alpha}(x)-\alpha_{0}(x)]^{2}F_{0}(dx)  &
=O_{p}\left(  \frac{(d_{K}+1)K}{n}+K^{-2\zeta_{\alpha}}\right)  .
\end{align*}

\bigskip

Proof: Note that for $\varepsilon_{i}=y_{i}-\gamma_{0}(x_{i})$, $E[\varepsilon
_{i}^{2}|x_{i}]=Var(y_{i}|x_{i})\leq C$. Note that $\hat{1}\hat{\Sigma}%
^{-2}\leq4I$ in the positive semi-definite semi-order so that%
\[
E[\hat{1}\left\Vert \hat{\Delta}_{1}\right\Vert ^{2}]\leq4E[\hat{h}_{1}%
^{T}\hat{h}_{1}]=\frac{4}{\hat{n}^{2}}\sum_{i,j\in\hat{I}}E[p(x_{i}%
)^{T}p(x_{j})\varepsilon_{i}\varepsilon_{j}]=\frac{4}{\hat{n}}E[\left\Vert
p(x_{i})\right\Vert ^{2}\varepsilon_{i}^{2}]\leq\frac{4C}{\hat{n}}E[\left\Vert
p(x_{i})\right\Vert ^{2}]=O(\frac{K}{n}).
\]
The first conclusion then follows by the Markov inequality. Next, we have
$\sup_{x}\left\vert \gamma_{K}(x)-\gamma_{0}(x_{i})\right\vert =O(K^{-\zeta})$
and hence for
\[
E[\hat{1}\left\Vert \hat{\Delta}_{2}\right\Vert ^{2}]\leq4E[\hat{h}_{2}%
^{T}\hat{h}_{2}]=\frac{4}{\hat{n}^{2}}\sum_{i,j\in\hat{I}}E[p(x_{i}%
)^{T}p(x_{j})r_{i}r_{j}]=\frac{4}{\hat{n}}E[\left\Vert p(x_{i})\right\Vert
^{2}]O(K^{-2\zeta_{\gamma}})=O\left(  K^{-2\zeta_{\gamma}}\frac{K}{n}\right)
,
\]
so the second equality also follows by the Markov inequality. Next, note that
\[
E[\eta_{i}^{T}\eta_{i}]\leq2E[v(z_{i})^{T}v(z_{i})]+2E[\alpha_{0}(x_{i}%
)^{2}\left\Vert p(x_{i})\right\Vert ^{2}]=O(K(d_{K}+1)).
\]
Then we have%
\[
E[\tilde{1}\left\Vert \tilde{\Delta}_{1}^{\alpha}\right\Vert ^{2}%
]\leq4E[\tilde{h}_{1}^{\alpha T}\tilde{h}_{1}^{\alpha}]=\frac{4}{\tilde{n}%
^{2}}\sum_{i,j\in\tilde{I}}E[\eta_{i}^{T}\eta_{j}]=\frac{4}{\hat{n}}E[\eta
_{i}^{T}\eta_{i}]=O(\frac{K(d_{K}+1)}{n}),
\]
so the third conclusion follows from the Markov inequality. The fourth
conclusion follows exactly like the second conclusion. the fifth and sixth
conclusions follow by eq. (\ref{coeff}) and the triangle inequality.

Next, note that by independence of the observations%
\begin{align*}
E[\hat{1}\hat{\Delta}_{1}\hat{\Delta}_{1}^{T}|X,\hat{Z}^{c}]  &  =\hat{1}%
\hat{\Sigma}^{-1}E[\hat{h}_{1}\hat{h}_{1}^{T}|X]\hat{\Sigma}^{-1}=\hat{1}%
\hat{\Sigma}^{-1}\left\{  \frac{1}{\hat{n}^{2}}\sum_{i,j\in\hat{I}}%
p(x_{i})p(x_{j})^{T}E[\varepsilon_{i}\varepsilon_{j}|X]\right\}  \hat{\Sigma
}^{-1}\\
&  =\hat{1}\hat{\Sigma}^{-1}\left\{  \frac{1}{\hat{n}^{2}}\sum_{i\in\hat{I}%
}p(x_{i})p(x_{i})^{T}E[\varepsilon_{i}^{2}|x_{i}]\right\}  \hat{\Sigma}%
^{-1}\leq\hat{1}\frac{C}{\hat{n}}\hat{\Sigma}^{-1}\leq\frac{2C}{n}I,
\end{align*}
giving the seventh conclusion.

Next, note that $\int p(x)[\gamma_{K}(x)-\gamma_{0}(x)]F_{0}(dx)=0$, so that%
\begin{align*}
\hat{1}\int[\hat{\gamma}(x)-\gamma_{0}(x)]^{2}F_{0}(dx)  &  =\hat{1}\int%
[\hat{\gamma}(x)-\gamma_{K}(x)+\gamma_{K}(x)-\gamma_{0}(x)]^{2}F_{0}(dx)\\
&  =\hat{1}\left\Vert \hat{\delta}-\delta\right\Vert ^{2}+\hat{1}%
K^{-2\zeta_{\gamma}}=O_{p}\left(  \frac{K}{n}+K^{-2\zeta_{\gamma}}\right)  ,
\end{align*}
giving the eighth conclusion. The last conclusion follows similarly.
\textit{Q.E.D.}

\bigskip

Next, we give an important intermediate result:

\bigskip

\textsc{Lemma A3}: \textit{If Assumptions 1-6 are satisfied then}
\[
\hat{1}\int[m(z,\hat{\gamma})-m(z,\gamma_{0})]^{2}F_{0}(dz)=O_{p}\left(
\frac{d_{K}K}{n}+K^{-2\zeta_{m}}\right)  .
\]
Proof: By linearity of $m(z,\gamma)-m(z,0)$, we have $m(z,\hat{\gamma
})-m(z,\gamma_{K})=v(z)^{T}(\hat{\delta}-\delta).$ Then by Lemma A2,%
\begin{align*}
\hat{1}\int[m(z,\hat{\gamma})-m(z,\gamma_{0})]^{2}F_{0}(dz)  &  \leq2\hat
{1}\int[m(z,\hat{\gamma})-m(z,\gamma_{K})]^{2}F_{0}(dz)+2\hat{1}%
\int[m(z,\gamma_{K})-m(z,\gamma_{0})]^{2}F_{0}(dz)\\
&  \leq2\hat{1}(\hat{\delta}-\delta)^{T}E[v(z_{i})v(z_{i})^{T}](\hat{\delta
}-\delta)+O(K^{-2\zeta_{m}})\\
&  \leq2d_{K}\hat{1}\left\Vert \hat{\delta}-\delta\right\Vert ^{2}%
+O(K^{-2\zeta_{m}})=O_{p}\left(  \frac{d_{K}K}{n}+K^{-2\zeta_{m}}\right)
.\text{ }Q.E.D.
\end{align*}

\bigskip

The proof of the results for the doubly robust estimators will make use of a
few Lemmas, that we now state.

\bigskip

\textsc{Lemma A4:}\textit{ If Assumptions 1-6 are satisfied then the
hypotheses of Lemma 6 are satisfied with }%
\[
\Delta_{n}^{m}=\sqrt{\frac{d_{K}K}{n}}+K^{-\zeta_{m}},\Delta_{n}^{\gamma
}=\sqrt{\frac{K}{n}}+K^{-\zeta_{\gamma}},\Delta_{n}^{\alpha}=\sqrt{\frac
{d_{K}K}{n}}+K^{-\zeta_{\alpha}}.
\]

\bigskip

Proof: The first conclusion follows by Lemma A3 and the second and third by
parts viii) and ix) of Lemma A2. $Q.E.D$.

\bigskip

\textsc{Lemma A5:}\textit{ If Assumptions 1-6 are satisfied and }$\hat{\gamma
}_{\ell}$ and $\tilde{\alpha}_{\ell}$ are computed from distinct samples then
for $\bar{\Sigma}=\bar{F}\{p(x)p(x)^{T}\}$%
\[
\sqrt{n}\bar{F}\{[\tilde{\alpha}_{\ell}(x)-\alpha_{0}(x)][\hat{\gamma}_{\ell
}(x)-\gamma_{0}(x)]\}=\sqrt{n}\hat{\Delta}_{2}^{T}\bar{\Sigma}\tilde{\Delta
}_{1}^{\alpha}+O_{p}(\bar{\Delta}_{n}^{\ast}+\Delta_{n}^{m}).
\]

\bigskip

Proof: Let $\bar{h}_{2}=\bar{F}\{p(x)[\gamma_{K}(x)-\gamma_{0}(x)]\}$ and
$\bar{h}_{2}^{\alpha}=\bar{F}\{p(x)[\alpha_{K}(x)-\alpha_{0}(x)]\}.$ Note that%
\begin{align*}
&  \bar{F}\{[\tilde{\alpha}_{\ell}(x)-\alpha_{0}(x)][\hat{\gamma}_{\ell
}(x)-\gamma_{0}(x)]\}\\
&  =\bar{F}\{[p(x)^{T}(\tilde{\delta}_{\alpha}-\delta_{\alpha})+\alpha
_{K}(x)-\alpha_{0}(x)][p(x)^{T}(\hat{\delta}-\delta)+\gamma_{K}(x)-\gamma
_{0}(x)]\}\\
&  =(\hat{\delta}-\delta)^{T}\bar{\Sigma}(\tilde{\delta}_{\alpha}%
-\delta_{\alpha})+(\hat{\delta}-\delta)^{T}\bar{h}_{2}^{\alpha}+(\hat{\delta
}_{\alpha}-\delta_{\alpha})^{T}\bar{h}_{2}+\bar{F}\{[\alpha_{K}(x)-\alpha
_{0}(x)][\gamma_{K}(x)-\gamma_{0}(x)]\}.
\end{align*}
By the Markov inequality%
\begin{equation}
\sqrt{n}\bar{F}\{[\alpha_{K}(x)-\alpha_{0}(x)][\gamma_{K}(x)-\gamma
_{0}(x)]\}=O_{p}(\sqrt{n}K^{-\zeta_{\gamma}-\zeta_{\alpha}}). \label{termbias}%
\end{equation}
Note that
\[
E[\tilde{h}_{2}^{\alpha}(\tilde{h}_{2}^{\alpha})^{T}]=\frac{1}{\bar{n}%
}E[p(x_{i})p(x_{i})^{T}(r_{i}^{\alpha})^{2}]\leq C\frac{1}{n}I.
\]
Therefore by Lemma A2 we have%
\[
E[\{\hat{1}(\hat{\delta}-\delta)^{T}\bar{h}_{2}^{\alpha}\}^{2}|\bar{Z}%
^{c}]=\hat{1}(\hat{\delta}-\delta)^{T}E[\tilde{h}_{2}^{\alpha}(\tilde{h}%
_{2}^{\alpha})^{T}](\hat{\delta}-\delta)\leq C\hat{1}\frac{1}{n}\left\Vert
\hat{\delta}-\delta\right\Vert ^{2}=O_{p}(\frac{K}{n^{2}}).
\]
Then by the Markov inequality it follows that%
\begin{equation}
\sqrt{n}(\hat{\delta}-\delta)^{T}\bar{h}_{2}^{\alpha}=O_{p}(\sqrt{\frac{K}{n}%
}). \label{term var1}%
\end{equation}
It follows similarly that%
\begin{equation}
\sqrt{n}(\hat{\delta}_{\alpha}-\delta_{\alpha})^{T}\bar{h}_{2}=O_{p}%
(\sqrt{\frac{d_{K}K}{n}}). \label{term var2}%
\end{equation}

Next, note that
\[
(\hat{\delta}-\delta)^{T}\bar{\Sigma}(\tilde{\delta}_{\alpha}-\delta_{\alpha
})=\hat{\Delta}_{1}^{T}\bar{\Sigma}(\tilde{\delta}_{\alpha}-\delta_{\alpha
})+\hat{\Delta}_{2}^{T}\bar{\Sigma}\tilde{\Delta}_{2}^{\alpha}+\hat{\Delta
}_{2}^{T}\bar{\Sigma}\tilde{\Delta}_{1}^{\alpha}.
\]
Let $\bar{1}$ be the event that $\lambda_{\max}(\bar{\Sigma})\leq2$. Then by
conclusion vii) of Lemma A2, and $\bar{1}$, $\hat{1}$, and $\tilde{1}$ all
functions of $X$ we have%
\begin{align*}
E[\bar{1}\hat{1}\tilde{1}\{\hat{\Delta}_{1}^{T}\bar{\Sigma}(\tilde{\delta
}_{\alpha}-\delta_{\alpha})\}^{2}|X,\hat{Z}^{c}]  &  =\bar{1}\hat{1}\tilde
{1}(\tilde{\delta}_{\alpha}-\delta_{\alpha})^{T}\bar{\Sigma}E[\hat{\Delta}%
_{1}\hat{\Delta}_{1}^{T}|X,\hat{Z}^{c}]\bar{\Sigma}(\tilde{\delta}_{\alpha
}-\delta_{\alpha})\\
&  \leq C\frac{1}{n}\bar{1}\tilde{1}(\tilde{\delta}_{\alpha}-\delta_{\alpha
})^{T}\bar{\Sigma}^{2}(\tilde{\delta}_{\alpha}-\delta_{\alpha})\leq4C\frac
{1}{n}\tilde{1}\left\Vert \tilde{\delta}_{\alpha}-\delta_{\alpha}\right\Vert
^{2}\\
&  =O_{p}(\frac{d_{K}K}{n^{2}}).
\end{align*}
Therefore we have%
\begin{equation}
\sqrt{n}\hat{\Delta}_{1}^{T}\bar{\Sigma}(\tilde{\delta}_{\alpha}%
-\delta_{\alpha})=O_{p}(\sqrt{\frac{d_{K}K}{n}}). \label{term var3}%
\end{equation}
Finally, note that by the Cauchy-Schwartz inequality%
\[
\hat{1}\hat{\Delta}_{2}^{T}\hat{\Delta}_{2}\leq\hat{1}2\hat{h}_{2}^{T}%
\hat{\Sigma}^{-1}\hat{h}_{2}\leq2\hat{F}\{[\gamma_{K}(x)-\gamma_{0}%
(x)]^{2}\}=O_{p}(K^{-2\zeta_{\gamma}}).
\]
It follows similarly that $\hat{1}(\hat{\Delta}_{2}^{\alpha})^{T}(\hat{\Delta
}_{2}^{\alpha})=O_{p}(K^{-2\zeta_{\alpha}})$ so that%
\begin{equation}
\sqrt{n}\bar{1}\hat{1}\tilde{1}\hat{\Delta}_{2}^{T}\bar{\Sigma}\tilde{\Delta
}_{2}^{\alpha}\leq2\sqrt{n}\sqrt{\hat{1}\hat{\Delta}_{2}^{T}\hat{\Delta}_{2}%
}\sqrt{\tilde{1}(\tilde{\Delta}_{2}^{\alpha})^{T}(\tilde{\Delta}_{2}^{\alpha
})}=O_{p}(\sqrt{n}K^{-\zeta_{\gamma}-\zeta_{\alpha}}). \label{term bias2}%
\end{equation}
The conclusion then follows by eqs. (\ref{termbias}), (\ref{term var1}),
(\ref{term var2}), (\ref{term var3}), (\ref{term bias2}), and the triangle
inequality. Q.E.D.

\bigskip

\textbf{Proof of Theorem 2: }It follows by Lemma A2 that the first hypothesis
of Lemma 1 is satisfied with $\Delta_{n}^{m}=\sqrt{d_{K}/n}+K^{-\zeta_{m}}\,$.
Let%
\[
\bar{m}(\gamma)=\int[m(z,\gamma)-m(z,\gamma)]F_{0}(dz)=E[\alpha_{0}%
(x_{i})\gamma(x_{i})],
\]
where the first equality is a definition and the second follows by Assumption
1. Then the first conclusion of Lemma 1 holds.

Next let $n=\hat{n}_{\ell}$ and $\hat{\gamma}=\hat{\gamma}_{\ell}$ for some
$\ell$ and $\phi(z)=\alpha_{0}(x)[y-\gamma_{0}(x)]$. Then it follows as in
Ichimura and Newey (2017), pp. 29 that%
\begin{align}
\hat{1}\sqrt{n}[\bar{m}(\hat{\gamma})-\beta_{0}-\frac{1}{n}\sum_{i=1}^{n}%
\phi(z_{i})]  &  =\hat{1}\left(  \hat{R}_{1}+\hat{R}_{2}+\hat{R}_{3}\right)
,\hat{R}_{1}=\sqrt{n}E[\alpha_{0}(x_{i})\{\gamma_{K}(x_{i})-\gamma_{0}%
(x_{i})\}],\label{plugexp}\\
\hat{R}_{2}  &  =\sqrt{n}v^{T}\hat{\Sigma}^{-1}\hat{h}_{2},\hat{R}_{3}%
=\frac{1}{\sqrt{n}}\sum_{i=1}^{n}[\alpha_{K}(x_{i})-\alpha_{0}(x_{i}%
)][y_{i}-\gamma_{0}(x_{i})].\nonumber
\end{align}
By $\gamma_{K}(x_{i})-\gamma_{0}(x_{i})$ orthogonal to $p(x_{i})$ in the
population and the Cauchy-Schwartz inequality,%
\begin{align*}
\left\vert \hat{R}_{1}\right\vert  &  =\sqrt{n}\left\vert E[\{\alpha_{0}%
(x_{i})-\alpha_{K}(x_{i})\}\{\gamma_{0}(x_{i})-\gamma_{K}(x_{i})\}]\right\vert
\leq\sqrt{n}\{E[\{\alpha_{0}(x_{i})-\alpha_{K}(x_{i})\}^{2}]E[\{\gamma
_{0}(x_{i})-\gamma_{K}(x_{i})\}^{2}]\}^{1/2}\\
&  =O(\sqrt{n}K^{-\zeta_{\gamma}-\zeta_{\alpha}})=O(\bar{\Delta}_{n}^{\ast}).
\end{align*}
Also,
\[
E[\hat{R}_{3}^{2}]=E[\{\alpha_{K}(x_{i})-\alpha(x_{i})\}^{2}\varepsilon
_{i}^{2}]\leq CE[\{\alpha_{K}(x_{i})-\alpha(x_{i})\}^{2}]=O(K^{-2\zeta
_{\alpha}}),
\]
so by the Markov inequality,%
\[
\hat{R}_{3}=O_{p}(K^{-\zeta_{\alpha}})=O_{p}(\bar{\Delta}_{n}^{\ast}).
\]

Next, note that $\hat{R}_{2}=\hat{R}_{21}+\hat{R}_{22}$ where $\hat{R}%
_{21}=v^{T}\hat{h}_{2}$ and $\hat{R}_{22}=\sqrt{n}v^{T}(\hat{\Sigma}%
^{-1}-I)\hat{h}_{2}.$ As noted following Assumption 4, $\sup_{x}|\gamma
_{K}(x)-\gamma_{0}(x)|=O(K^{-\zeta_{\gamma}})$, so that
\[
E[\hat{R}_{21}^{2}]=v^{T}E[p(x_{i})p(x_{i})^{T}r_{i}^{2}]v\leq O(K^{-2\zeta
_{\gamma}})v^{T}v\leq O(K^{-2\zeta_{\gamma}})E[\alpha_{0}(x_{i})^{2}%
]=O(K^{-2\zeta_{\gamma}}).
\]
Then by the Markov inequality%
\[
\hat{R}_{21}=O_{p}(K^{-\zeta_{\gamma}})=O_{p}(\bar{\Delta}_{n}^{\ast}).
\]
Finally, note that $(\hat{\Sigma}^{-1}-I)\hat{h}_{2}=\hat{U}+\hat{W}$ for
$\hat{U}$ and $\hat{W}$ defined in the statement of Lemma A1, so that for any
$\Delta>0$ we have
\begin{align*}
\hat{1}\hat{R}_{22}^{2}  &  =\hat{1}n\cdot v^{T}(\hat{U}+\hat{W})(\hat{U}%
+\hat{W})v\leq2\hat{1}n\cdot v^{T}(\hat{U}\hat{U}^{T}+\hat{W}\hat{W})v\\
&  \leq CK^{-2\zeta_{\gamma}}[\ln(n)]^{2}+O_{p}(n^{-\Delta+1})\text{,}%
\end{align*}
for any $C$. It then follows by eq. (\ref{plugexp}) and the triangle
inequality that%
\[
\sqrt{n}[\bar{m}(\hat{\gamma})-\beta_{0}]=\frac{1}{\sqrt{n}}\sum_{i=1}^{n}%
\phi(z_{i})+O_{p}(\bar{\Delta}_{n}^{\ast}+K^{-\zeta_{\gamma}}\ln(n)).
\]
The first conclusion then follows from the second conclusion of Lemma 1. The
second conclusion follows by $\Delta_{n}^{m}=C\sqrt{K/n}+K^{-\zeta_{\gamma}%
}=O(\bar{\Delta}_{n}^{\ast})$ when $d_{K}$ is bounded and $\zeta_{m}%
=\zeta_{\gamma}$. $Q.E.D.$

\bigskip

\textbf{Proof of Corollary 3: }To prove this result it suffices to show that
Assumptions 4 and 6 are satisfied in each of the examples with $d_{K}$ bounded
and $\zeta_{m}=\zeta_{\gamma}.$

For the conditional covariance $\alpha_{0}(x)=-E[a_{i}|x_{i}=x]$. This being
Holder of order $s_{\alpha}$ is a hypothesis. Also, $v(z)=a\cdot p(x),$ so
that
\[
E[v(z_{i})v(z_{i})^{T}]=E[a_{i}^{2}p(x_{i})p(x_{i})^{T}]\leq CE[p(x_{i}%
)p(x_{i})^{T}]\leq CI
\]
by $E[a_{i}^{2}|x_{i}]$ bounded. Also $\zeta_{m}=\zeta_{\gamma}$ by%
\[
E[\left\{  m(z_{i},\gamma_{K})-m(z_{i},\gamma_{0})\right\}  ^{2}]=E[a_{i}%
^{2}\{\gamma_{K}(x_{i})-\gamma_{0}(x_{i})\}^{2}]\leq CK^{-2\zeta_{\gamma}}.
\]

For the missing data mean $\alpha_{0}(x)=a/\pi_{0}(w)$ is Holder of order
$s_{\alpha}$ by $\pi_{0}(w_{i})$ being bounded away from zero and Holder of
order $s_{\alpha}.$ Furthermore $v(z)=q(w)$, so that by Assumption 3,
\[
E[v(z_{i})v(z_{i})^{T}]=E[q(w_{i})q(w_{i})^{T}]\leq CI,
\]
and by $a_{i}$ bounded and $\pi_{0}(w_{i})$ bounded away from zero,%
\begin{align*}
E[\left\{  m(z_{i},\gamma_{K})-m(z_{i},\gamma_{0})\right\}  ^{2}]  &
=E[\{q(w_{i})^{T}\delta-E[y_{i}|a_{i}=1,w_{i}]\}^{2}]\\
&  =E[\frac{a_{i}}{\pi_{0}(w_{i})}\{\gamma_{K}(x_{i})-\gamma_{0}(x_{i}%
)\}^{2}]\\
&  \leq CE[\{\gamma_{K}(x_{i})-\gamma_{0}(x_{i})\}^{2}]\leq CK^{-2\zeta
_{\gamma}}.
\end{align*}

For the average derivative example $\alpha_{0}(x)=\omega(x)/f_{0}(x)$ which is
Holder of order $s_{\alpha}$ by each of $\omega(x)$ and $f_{0}(x)$ being
Holder of order $s_{\alpha}$ and by $f_{0}(x)$ bounded away from zero where
$\omega(x)$ is non zero. Furthermore $v(z)=\int\omega(x)p(x)dx,$ so that by
Cauchy-Schwartz,
\begin{align*}
E[v(z_{i})v(z_{i})^{T}]  &  =\int\omega(x)p(x)dx\int\omega(x)p(x)^{T}dx\\
&  =E[\alpha_{0}(x_{i})p(x_{i})]E[\alpha_{0}(x_{i})p(x_{i})^{T}]\\
&  \leq E[\alpha_{0}(x_{i})^{2}]E[p(x_{i})p(x_{i})^{T}]\leq CI.
\end{align*}
Furthermore,
\begin{align*}
E[\left\{  m(z_{i},\gamma_{K})-m(z_{i},\gamma_{0})\right\}  ^{2}]  &
=\{\int\omega(x)[\gamma_{K}(x)-\gamma_{0}(x)]dx\}^{2}\\
&  =E[\alpha_{0}(x_{i})\{\gamma_{K}(x_{i})-\gamma_{0}(x_{i})\}]^{2}\\
&  \leq E[\alpha_{0}(x_{i})^{2}]E[\{\gamma_{K}(x_{i})-\gamma_{0}(x_{i}%
)\}^{2}]=O(K^{-2\zeta_{\gamma}}).\text{ }Q.E.D.
\end{align*}

\bigskip

\textbf{Proof of Theorem 4:} The conclusion follows from Lemma 1 and Theorem 8
of Ichimura and Newey (2017) similarly to the proof of Theorem 2 above, with
the conclusion of Theorem 8 of Ichimura and Newey (2017) replacing the
argument following eq. (\ref{plugexp}) in the proof of Theorem 2. $Q.E.D.$

\bigskip

\textbf{Proof of Theorem 5:} Let $\hat{\lambda}_{\ell}(x)$ denote the series
regression of $u_{i}=y-a_{i}^{T}\beta_{0}$ on $p(x_{i})$ in the $\hat{I}%
_{\ell}$ sample. By a standard formula for instrumental variables estimation
and series estimation,
\begin{align}
\sqrt{n}(\hat{\beta}-\beta_{0})  &  =\hat{H}^{-1}\frac{1}{\sqrt{n}}\sum
_{\ell=1}^{L}\sum_{i\in I_{\ell}}[a_{i}-\hat{\alpha}_{\ell}(x_{i})]\left\{
y_{i}-\hat{\gamma}_{\ell}(x_{i})-[a_{i}-\hat{\alpha}_{\ell}(x_{i})]^{T}%
\beta_{0}\right\} \label{plex}\\
&  =\hat{H}^{-1}\frac{1}{\sqrt{n}}\sum_{\ell=1}^{L}\sum_{i\in I_{\ell}}%
[a_{i}-\hat{\alpha}_{\ell}(x_{i})]\left[  u_{i}-\hat{\lambda}_{\ell}%
(x_{i})\right] \nonumber
\end{align}
Assume for the moment that $a_{i}$ is a scalar and let $y_{i}=u_{i}$. Then
$\sum_{\ell=1}^{L}\sum_{i\in I_{\ell}}[a_{i}-\hat{\alpha}_{\ell}%
(x_{i})]\left[  u_{i}-\hat{\lambda}_{\ell}(x_{i})\right]  /n$ is the doubly
robust estimator with $m(z,\gamma)=a[y-\gamma(x)],$ i.e. for the expected
conditional covariance. It then follows as in the proof of Corollary 3 that
$\max\{\Delta_{n}^{m},\Delta_{n}^{\gamma},\Delta_{n}^{\alpha}\}\leq
C\bar{\Delta}_{n}^{\ast}.$ Then by Lemmas 6 and A5, for $\varphi
(z)=[a_{i}-\alpha_{0}(x_{i})]\varepsilon_{i}$,
\[
\frac{1}{\sqrt{n}}\sum_{\ell=1}^{L}\sum_{i\in I_{\ell}}[a_{i}-\hat{\alpha
}_{\ell}(x_{i})]\left[  u_{i}-\hat{\lambda}_{\ell}(x_{i})\right]  =\frac
{1}{\sqrt{n}}\sum_{i=1}^{n}\varphi(z_{i})+O_{p}(\bar{\Delta}_{n}^{\ast}%
)+\sqrt{n}\hat{\Delta}_{2}^{T}\bar{\Sigma}\tilde{\Delta}_{1}^{\alpha}.
\]
Note that here $\alpha_{0}(x_{i})=-E[a_{i}|x_{i}]$ so that
\[
\tilde{h}_{1}^{\alpha}=\tilde{F}\{v(z)-\alpha_{0}(x)p(x)\}=-\tilde
{F}\{[a-\alpha_{0}(x)]p(x)\}.
\]
Then we have
\[
E[\tilde{1}\tilde{\Delta}_{1}^{\alpha}\tilde{\Delta}_{1}^{\alpha T}%
|X,\tilde{Z}^{c}]=\tilde{1}\frac{1}{\tilde{n}}\tilde{\Sigma}^{-1}\tilde
{F}\{p(x)p(x)^{T}Var(a_{i}|x_{i}=x)\}\tilde{\Sigma}^{-1}\leq\frac{C}{n}I.
\]
Therefore it follows by Lemma A2 that%
\[
E[\hat{1}\tilde{1}(\hat{\Delta}_{2}^{T}\bar{\Sigma}\tilde{\Delta}_{1}^{\alpha
})^{2}|X,\tilde{Z}^{c}]=\hat{1}\hat{\Delta}_{2}^{T}E[\tilde{1}\tilde{\Delta
}_{1}^{\alpha}\tilde{\Delta}_{1}^{\alpha T}|X,\tilde{Z}^{c}]\hat{\Delta}%
_{2}\leq\hat{1}\frac{C}{n}\hat{\Delta}_{2}^{T}\hat{\Delta}_{2}=o_{p}%
(K/n^{2}).
\]
Then by the Markov inequality%
\[
\sqrt{n}\hat{\Delta}_{2}^{T}\bar{\Sigma}\tilde{\Delta}_{1}^{\alpha}%
=o_{p}\left(  \sqrt{\frac{K}{n}}\right)  =O_{p}(\bar{\Delta}_{n}^{\ast}).
\]
Consequently we have%
\begin{equation}
\frac{1}{\sqrt{n}}\sum_{\ell=1}^{L}\sum_{i\in I_{\ell}}[a_{i}-\hat{\alpha
}_{\ell}(x_{i})]\left[  u_{i}-\hat{\lambda}_{\ell}(x_{i})\right]  =\frac
{1}{\sqrt{n}}\sum_{i=1}^{n}\varphi(z_{i})+O_{p}(\bar{\Delta}_{n}^{\ast}).
\label{plscore}%
\end{equation}

Next, note that%
\begin{align*}
\bar{F}\{[a-\tilde{\alpha}_{\ell}(x)][a-\hat{\alpha}_{\ell}(x)]\}  &  =\bar
{F}\{[a-\alpha_{0}(x)+\alpha_{0}(x)-\tilde{\alpha}_{\ell}(x)][a-\alpha
_{0}(x)+\alpha_{0}(x)-\hat{\alpha}_{\ell}(x)]\}\\
&  =\bar{F}\{[a-\alpha_{0}(x)]^{2}+[a-\alpha_{0}(x)][\alpha_{0}(x)-\tilde
{\alpha}_{\ell}(x)]\\
&  +[a-\alpha_{0}(x)][\alpha_{0}(x)-\hat{\alpha}_{\ell}(x)]+[\alpha
_{0}(x)-\tilde{\alpha}_{\ell}(x)][\alpha_{0}(x)-\hat{\alpha}_{\ell}(x)]\}
\end{align*}
Note that by Lemma A3 and $E[a_{i}^{2}|x_{i}]$ bounded,
\begin{align*}
\tilde{1}E[(\bar{F}\{[a-\alpha_{0}(x)][\alpha_{0}(x)-\tilde{\alpha}_{\ell
}(x)]\})^{2}|\tilde{Z}]  &  =\frac{\tilde{1}}{\bar{n}}\int[a-\alpha
_{0}(x)]^{2}[\alpha_{0}(x)-\tilde{\alpha}(x)]^{2}F_{0}(dz)\\
&  \leq C\frac{\tilde{1}}{n}\int E[a_{i}^{2}|x_{i}=x][\alpha_{0}%
(x)-\tilde{\alpha}(x)]^{2}F_{0}(dz)\\
&  \leq C\frac{\tilde{1}}{n}\int[\alpha_{0}(x)-\tilde{\alpha}(x)]^{2}%
F_{0}(dz)=O_{p}\left(  \frac{1}{n}\left(  \frac{K}{n}+K^{-2\zeta_{\gamma}%
}\right)  \right)  .
\end{align*}
so that by the Markov inequality it follows that%
\begin{equation}
\bar{F}\{[a-\alpha_{0}(x)][\alpha_{0}(x)-\tilde{\alpha}_{\ell}(x)]\}=O_{p}%
(\bar{\Delta}_{n}^{\ast}). \label{Term1}%
\end{equation}
\qquad\ It follows similarly that%
\begin{equation}
\bar{F}\{[a-\alpha_{0}(x)][\alpha_{0}(x)-\hat{\alpha}_{\ell}(x)]\}=O_{p}%
(\bar{\Delta}_{n}^{\ast}). \label{term2}%
\end{equation}
Also, by the Cauchy-Schwartz inequality
\[
\hat{1}\tilde{1}\left\vert \bar{F}\{[\alpha_{0}(x)-\tilde{\alpha}_{\ell
}(x)][\alpha_{0}(x)-\hat{\alpha}_{\ell}(x)]\}\right\vert \leq(\tilde{1}\bar
{F}\{[\alpha_{0}(x)-\tilde{\alpha}_{\ell}(x)]^{2}\})^{1/2}(\hat{1}\bar
{F}\{[\alpha_{0}(x)-\hat{\alpha}_{\ell}(x)]^{2}\})^{1/2}.
\]
Also,%
\[
E[\tilde{1}\bar{F}\{[\alpha_{0}(x)-\tilde{\alpha}_{\ell}(x)]^{2}\}|\tilde
{Z}]=\tilde{1}\int[\tilde{\alpha}_{\ell}(x)-\alpha_{0}(x)]^{2}F_{0}%
(dx)=O_{p}\left(  \frac{K}{n}+K^{-2\zeta_{\gamma}}\right)  ,
\]
so that $\tilde{1}\bar{F}\{[\alpha_{0}(x)-\tilde{\alpha}_{\ell}(x)]^{2}%
\}=O_{p}(K/n+K^{-2\zeta_{\gamma}}).$ It follows similarly that $\hat{1}\bar
{F}\{[\alpha_{0}(x)-\hat{\alpha}_{\ell}(x)]^{2}\}=O_{p}(K/n+K^{-2\zeta
_{\gamma}}),$ so that%
\begin{equation}
=\bar{F}\{[\alpha_{0}(x)-\tilde{\alpha}_{\ell}(x)][\alpha_{0}(x)-\hat{\alpha
}_{\ell}(x)]\}=O_{p}(\bar{\Delta}_{n}^{\ast}) \label{Term3}%
\end{equation}
Also, note that by $E[\left\Vert a_{i}\right\Vert ^{4}]<\infty,$%
\[
\bar{F}\{[a-\alpha_{0}(x)]^{2}\}=E[\{a-\alpha_{0}(x)\}^{2}]+O_{p}\left(
\frac{1}{\sqrt{n}}\right)  =E[\{a_{i}-\alpha_{0}(x_{i})\}^{2}]+O_{p}%
(\bar{\Delta}_{n}^{\ast}).
\]
It then follows by eqs. (\ref{Term1}), (\ref{term2}), (\ref{Term3}) and the
triangle inequality that%
\[
\bar{F}\{[a-\tilde{\alpha}_{\ell}(x)][a-\hat{\alpha}_{\ell}(x)]\}=E[\{a_{i}%
-\alpha_{0}(x_{i})\}^{2}]+O_{p}(\bar{\Delta}_{n}^{\ast}).
\]
Applying this argument to each element of $\hat{H}=\sum_{\ell=1}^{L}\sum_{i\in
I_{\ell}}[a_{i}-\tilde{\alpha}_{\ell}(x_{i})][a_{i}-\hat{\alpha}_{\ell}%
(x_{i})]^{T}/n$ and each group of observations $I_{\ell}$ and summing up gives
$\hat{H}=H+O_{p}(\bar{\Delta}_{n}^{\ast}).$ It then follows by a standard
argument and nonsingularity of $H$ that
\begin{equation}
\hat{H}^{-1}=H^{-1}+O_{p}(\bar{\Delta}_{n}^{\ast}). \label{plhess}%
\end{equation}

Finally, it follows from eqs. (\ref{plex}), (\ref{plscore}), (\ref{plhess})
and from $\sum_{i=1}^{n}\varphi(z_{i})/\sqrt{n}=O_{p}(1)$ that%
\[
\sqrt{n}(\hat{\beta}-\beta_{0})=[H^{-1}+O_{p}(\bar{\Delta}_{n}^{\ast}%
)][\frac{1}{\sqrt{n}}\sum_{i=1}^{n}\varphi(z_{i})+O_{p}(\bar{\Delta}_{n}%
^{\ast})]=H^{-1}\frac{1}{\sqrt{n}}\sum_{i=1}^{n}\varphi(z_{i})+O_{p}%
(\bar{\Delta}_{n}^{\ast}).\text{ }Q.E.D.
\]

\bigskip

\textbf{Proof of Theorem 7:} By Lemmas 6 and A5 it suffices to show that
$\bar{1}\hat{1}\tilde{1}\sqrt{n}\hat{\Delta}_{2}^{T}\bar{\Sigma}\tilde{\Delta
}_{1}^{\alpha}=O_{p}(\Delta_{n}^{m}).$ Note that%
\[
\hat{1}\hat{\Delta}_{2}=\hat{1}\hat{h}_{2}+\hat{1}\hat{U}+\hat{1}\hat{W}.
\]
By $E[\hat{h}_{2}\hat{h}_{2}^{T}]\leq Cn^{-1}K^{-2\zeta_{\gamma}}I$ \ and
Lemma A2 iii) we have
\begin{align*}
E[\left(  \bar{1}\hat{1}\tilde{1}\sqrt{n}\hat{h}_{2}^{T}\bar{\Sigma}%
\tilde{\Delta}_{1}^{\alpha}\right)  ^{2}|\hat{Z}^{c}]  &  =n\bar{1}\tilde
{1}\left(  \tilde{\Delta}_{1}^{\alpha}\right)  ^{T}\bar{\Sigma}E[\hat{1}%
\hat{h}_{2}\hat{h}_{2}^{T}]\bar{\Sigma}\tilde{\Delta}_{1}^{\alpha}\leq
n\bar{1}\tilde{1}\left(  \tilde{\Delta}_{1}^{\alpha}\right)  ^{T}\bar{\Sigma
}E[\hat{h}_{2}\hat{h}_{2}^{T}]\bar{\Sigma}\tilde{\Delta}_{1}^{\alpha}\\
&  \leq CK^{-2\zeta_{\gamma}}\bar{1}\tilde{1}\left(  \tilde{\Delta}%
_{1}^{\alpha}\right)  ^{T}\bar{\Sigma}^{2}\tilde{\Delta}_{1}^{\alpha}%
=O_{p}\left(  \frac{(1+d_{K})K^{1-2\zeta_{\gamma}}}{n}\right)  =O_{p}%
((\Delta_{n}^{m})^{2}).
\end{align*}
Also, by the first conclusion of Lemma A1 and by Lemma A2 iii),%
\begin{align*}
E[\left(  \bar{1}\hat{1}\tilde{1}\sqrt{n}\hat{U}^{T}\bar{\Sigma}\tilde{\Delta
}_{1}^{\alpha}\right)  ^{2}|\hat{Z}^{c}]  &  =n\bar{1}\tilde{1}\left(
\tilde{\Delta}_{1}^{\alpha}\right)  ^{T}\bar{\Sigma}E[\hat{1}\hat{U}\hat
{U}^{T}]\bar{\Sigma}\tilde{\Delta}_{1}^{\alpha}\leq n\bar{1}\tilde{1}\left(
\tilde{\Delta}_{1}^{\alpha}\right)  ^{T}\bar{\Sigma}E[\hat{U}\hat{U}^{T}%
]\bar{\Sigma}\tilde{\Delta}_{1}^{\alpha}\\
&  \leq CK^{-2\zeta_{\gamma}}\ln(n)^{2}\bar{1}\tilde{1}\left(  \tilde{\Delta
}_{1}^{\alpha}\right)  ^{T}\bar{\Sigma}^{2}\tilde{\Delta}_{1}^{\alpha}%
=O_{p}\left(  \frac{(1+d_{K})K^{1-2\zeta_{\gamma}}[\ln(n)]^{2}}{n}\right)
=O_{p}((\Delta_{n}^{m})^{2}).
\end{align*}
Also by the second conclusion of Lemma A1 and Lemma A2 iii), for $\Delta>0$
large enough,
\[
\bar{1}\hat{1}\tilde{1}\sqrt{n}\hat{\Delta}_{2}^{T}\bar{\Sigma}\tilde{\Delta
}_{1}^{\alpha}=O_{p}(n^{(1/2)-\Delta}\sqrt{(1+d_{K})/n}))=O_{p}(\Delta_{n}%
^{m}).
\]
The conclusion then follows by the Markov and triangle inequalities.
\textit{Q.E.D.}

\bigskip

\textbf{Proof of Theorem 8: }By Lemmas 6 and A5 it suffices to show that
$\bar{1}\hat{1}\tilde{1}\sqrt{n}\hat{\Delta}_{2}^{T}\bar{\Sigma}\tilde{\Delta
}_{1}^{\alpha}=O_{p}(\Delta_{n}^{m}+\tilde{\Delta}_{n}).$ Note that%
\begin{align*}
\bar{1}\hat{1}\tilde{1}\sqrt{n}\hat{\Delta}_{2}^{T}\bar{\Sigma}\tilde{\Delta
}_{1}^{\alpha}  &  =T_{1}+T_{2}+T_{3},T_{1}=\bar{1}\hat{1}\tilde{1}\sqrt
{n}\hat{h}_{2}^{T}\bar{\Sigma}\tilde{\Delta}_{1}^{\alpha},\\
T_{2}  &  =\bar{1}\hat{1}\tilde{1}\sqrt{n}\hat{\Delta}_{2}^{T}(I-\hat{\Sigma
})\bar{\Sigma}\tilde{h}_{1}^{\alpha},\\
T_{3}  &  =\bar{1}\hat{1}\tilde{1}\sqrt{n}\hat{\Delta}_{2}^{T}(I-\hat{\Sigma
})\bar{\Sigma}(I-\tilde{\Sigma})\tilde{\Delta}_{1}^{\alpha}.
\end{align*}
By Lemma A2 iii),
\[
E[T_{1}^{2}|\hat{Z}^{c}]\leq\bar{1}\tilde{1}n(\tilde{\Delta}_{1}^{\alpha}%
)^{T}\bar{\Sigma}E[\hat{h}_{2}\hat{h}_{2}^{T}]\bar{\Sigma}\tilde{\Delta}%
_{1}^{\alpha}\leq CK^{-2\zeta_{\gamma}}\tilde{1}(\tilde{\Delta}_{1}^{\alpha
})^{T}\tilde{\Delta}_{1}^{\alpha}=O_{p}(K^{-2\zeta_{\gamma}}\left(  \Delta
_{n}^{m}\right)  ^{2}),
\]
so by the Markov inequality, $T_{1}=O_{p}(\Delta_{n}^{m}).$ By Lemma A2 ii),
\begin{align*}
E[T_{2}^{2}|\tilde{Z}^{c}]  &  \leq\hat{1}\sqrt{n}\hat{\Delta}_{2}^{T}%
(I-\hat{\Sigma})E[\tilde{h}_{1}^{\alpha}\left(  \tilde{h}_{1}^{\alpha}\right)
^{T}](I-\hat{\Sigma})\hat{\Delta}_{2}\leq Cd_{K}\hat{\Delta}_{2}^{T}%
(I-\hat{\Sigma})^{2}\hat{\Delta}_{2}\\
&  =O_{p}((1+d_{K})\frac{K^{1-2\zeta_{\gamma}}}{n}\frac{K\ln(K)}{n}).
\end{align*}
Note that by the Markov inequality and $K\ln(K)/n\longrightarrow0$ it follows
that $T_{2}=O_{p}(\bar{\Delta}_{n}^{\ast}+\Delta_{n}^{m}).$ Finally, by the
Caucy-Schwartz inequality and Lemma A2,%
\[
T_{3}=O_{p}(\sqrt{\frac{K^{3}\ln(K)(1+d_{K})}{n^{3}}}K^{(1/2)-\zeta_{\gamma}%
})=O_{p}(\tilde{\Delta}_{n}).
\]
The conclusion then follows by the triangle inequality. \textit{Q.E.D.}

\section*{Acknowledgements}

We appreciate the hospitality of the Cowles Foundation where much of the work
for this paper was accomplished. We also appreciate the comments of M.
Cattaneo, X. Chen, M. Jansson and seminar participants at UCL.

\bigskip

\setlength{\parindent}{-.5cm} \setlength{\parskip}{.1cm}

\begin{center}
\textbf{REFERENCES}
\end{center}

\textsc{Athey, S., G. Imbens, and S. Wager} (2017): "Efficient Inference of
Average Treatment Effects in High Dimensions via Approximate Residual
Balancing," \textit{Journal of the Royal Statistical Society, Series B,} forthcoming.

\textsc{Ayyagari, R. }(2010): Applications of Influence Functions to
Semiparametric Regression Models, Ph.D. Thesis, Harvard School of Public
Health, Harvard University.

\textsc{Belloni, A., V. Chernozhukov, D. Chetverikov, K. Kato} (2015):
\textquotedblleft Some New Asymptotic Theory for Least Squares Series:
Pointwise and Uniform Results,\textquotedblright\ \textit{Journal of
Econometrics} 186, 345--366.

\textsc{Bickel, P.J.} (1982): "On Adaptive Estimation," \textit{Annals of
Statistics} 10, 647-671.

\textsc{Bickel, P. and Y. Ritov} (1988): "Estimating Integrated Squared
Density Derivatives: Sharp Best Order of Convergence Estimates,"
\textit{Sankhya: The Indian Journal of Statistics}, \textit{Series A} 50, 381--393.

\textsc{Blomquist, S. and M. Dahlberg} (1999): "Small Sample Properties of
LIML and Jackknife IV Estimators: Experiments with Weak Instruments,"
\textit{Journal of Applied Econometrics }14, 69--88.

\textsc{Cattaneo, M.D., and M. Farrell} (2013): "Optimal Convergence Rates,
Bahadur Representation, and Asymptotic Normality of Partitioning Estimators,"
\textit{Journal of Econometrics} 174, 127-143.

\textsc{Cattaneo, M.D., and M. Jansson }(2017): "Kernel-Based Semiparametric
Estimators: Small Bandwidth Asymptotics and Bootstrap Consistency,"
\textit{Econometrica}, forthcoming.

\textsc{Cattaneo, M.D., M. Jansson, and X. Ma (2017): "}Two-step Estimation
and Inference with Possibly Many Included Covariates," working paper, Michigan.

\textsc{Chernozhukov, V., J.C. Escanciano, H. Ichimura, W.K. Newey, J.M.
Robins }(2016): "Locally Robust Semiparametric Estimation," arXiv 1608.00033.

\textsc{Chernozhukov, V., D. Chetverikov, M. Demirer, E. Duflo, C. Hansen,
W.K. Newey, J.M. Robins }(2017): "Double/Debiased Machine Learning for
Treatment and Structural Parameters," \textit{Econometrics Journal}, forthcoming.

\textsc{Donald, S.G. and W.K. Newey }(1994): \textquotedblleft Series
Estimation of Semilinear Models," \textit{Journal of Multivariate Analysis}
50, 30-40.

\textsc{Firpo, S. and C. Rothe} (2016): "Semiparametric Two-Step Estimation
Using Doubly Robust Moment Conditions," working paper.

\textsc{Gine, E. and R. Nickl} (2008): "A Simple Adaptive Estimator of the
Integrated Square of a Density," \textit{Bernoulli} 14, 47--61.

\textsc{Hahn, J. (1998):} "On the Role of the Propensity Score in Efficient
Semiparametric Estimation of Average Treatment Effects," \textit{Econometrica}
66, 315-331.

\textsc{Hirano, K., G. Imbens, and G. Ridder} (2003): "Efficient Estimation of
Average Treatment Effects Using the Estimated Propensity Score,"
\textit{Econometrica} 71: 1161--1189.

\textsc{Hirschberg, D.A and S. Wager} (2017): "Balancing Out Regression Error:
Efficient Treatment Effect Estimation without Smooth Propensities," arXiv:1712.00038.

\textsc{Ichimura, H. and W.K. Newey} (2017): "The Influence Function of
Semiparametric Estimators," CEMMAP working paper CWP06/17.

\textsc{Imbens G., J. Angrist, A. Krueger} (1999): "Jackknife Instrumental
Variables Estimation," \textit{Journal of Applied Econometrics }14, 57-67.

\textsc{Kandasamy, K., A. Krishnamurthy, B. Poczos, L. Wasserman, J. Robins
}(2015) "Nonparametric von Mises Estimators for Entropies, Divergences and
Mutual Informations," \textit{Advances in Neural Information Processing
Systems} 28 (NIPS 2015).

\textsc{Laurent, B.} (1996): "Efficient Estimation of Integral Functionals of
a Density," \textit{Annals of Statistics} 24, 659-681.

\textsc{Mukherjee, R., W.K. Newey, J.M. Robins} (2017): "Semiparametric
Efficient Empirical Higher Order Influence Function Estimators," arXiv:1705.07577.

\textsc{Newey, W.K.} (1994): "The Asymptotic Variance of Semiparametric
Estimators," \textit{Econometrica} 62, 1349-1382.

\textsc{Newey, W.K. }(1997): \textquotedblleft Convergence Rates and
Asymptotic Normality for Series Estimators,\textquotedblright\ \textit{Journal
of Econometrics }79, 147-168.

\textsc{Newey, W.K., F. Hsieh, {\small and} J.M. Robins} (1998):
\textquotedblleft Undersmoothing and Bias Corrected Functional Estimation,"
MIT Dept. of Economics working paper\ 72, 947-962.

\textsc{Newey, W.K., F. Hsieh, {\small and} J.M. Robins} (2004):
\textquotedblleft Twicing Kernels and a Small Bias Property of Semiparametric
Estimators,\textquotedblright\ \textit{Econometrica} 72, 947-962.

\textsc{Powell, J.L., J.H. Stock, and T.M. Stoker }(1989): "Semiparametric
Estimation of Index Coefficients," \textit{Econometrica} 57, 1403-1430.

\textsc{Robins, J.M. and A. Rotnitzky (1995): }"Semiparametric Efficiency in
Multivariate Regression Models with Missing Data," \textit{Journal of the
American Statististical Association} 90, 122--129.

\textsc{Robins, J.M., A. Rotnitzky, and M. van der Laan} \ (2000): "Comment on
'On Profile Likelihood'\ by S. A. Murphy and A. W. van der Vaart,"
\textit{Journal of the American Statistical Association} 95, 431-435.

\textsc{Robins, J., M. Sued, Q. Lei-Gomez, and A. Rotnitzky} (2007): "Comment:
Performance of Double-Robust Estimators When Inverse Probability' Weights Are
Highly Variable," \textit{Statistical Science} 22, 544--559.

\textsc{Robins, J.M., E.T. Tchetgen, L. Li, A. van der Vaart (2009):}
"Semiparametric Minimax Rates," \textit{Electronic Journal of Statistics }3, 1305--1321.

\textsc{Robins, J.M., L. Li, E. Tchetgen, and A. van der Vaart} (2008) "Higher
Order Influence Functions and Minimax Estimation of Nonlinear Functionals," in
\textit{IMS Collections Vol. 2, Probability and Statistics: Essays in Honor of
David A. Freedman, }D. Nolan and T. Speed (eds.), Beachwood, Ohio: Institute
of Mathematical Statistics, 335-421.

\textsc{Robins, J.M, P. Zhang, R. Ayyagari, R. Logan, E. Tchetgen, L. Li, T.
Lumley, A. van der Vaart A, HEI Health Review Committee }(2013): "New
Statistical Approaches to Semiparametric Regression with Application to Air
Pollution Research," Research Report Health Eff Instm 175:3-129.

\textsc{Robins, J.M., L. Li, R. Mukherjee, E. Tchetgen, A. van der Vaart}
(2017): "Minimax Estimation of a Functional on a Structured High Dimensional
Model," \textit{Annals of Statistics,} forthcoming.

\textsc{Rotnitzky, A. and J.M. Robins} (1995): "Semi-parametric Estimation of
Models for Means and Covariances in the Presence of Missing Data,"
\textit{Scandinavian Journal of Statistics} 22, 323--333.

\textsc{Rudelson, M.} (1999): "Random Vectors in the Isotropic Position,"
\textit{Journal of Functional Analysis} 164, 60-72.

\textsc{Scharfstein D.O., A. Rotnitzky, and J.M. Robins (1999): }Rejoinder to
\textquotedblleft Adjusting For Nonignorable Drop-out Using Semiparametric
Non-response Models,\textquotedblright\ \textit{Journal of the American
Statistical Association }94, 1135-1146.

\textsc{Stoker, T. }(1986): "Consistent Estimation of Scaled Coefficients,"
\textit{Econometrica} 54, 1461-1482.

\setlength{\parindent}{.0cm} \setlength{\parskip}{.1cm}

\end{document}